\documentclass{amsart}
\usepackage{graphicx}
\usepackage[all]{xy}
\usepackage{latexsym}
\usepackage{amssymb}
\usepackage{amsmath}
\usepackage{color}

\xyoption{arc}

 \newfam\cyrfam

  \font\tencyr=wncyr10

  \font\sevencyr=wncyr7

  \font\fivecyr=wncyr5

  \textfont\cyrfam=\tencyr \scriptfont\cyrfam=\sevencyr

    \scriptscriptfont\cyrfam=\fivecyr

  \newfam\cyifam

  \font\tencyi=wncyi10

  \font\sevencyi=wncyi7

  \font\fivecyi=wncyi5

  \textfont\cyifam=\tencyi \scriptfont\cyifam=\sevencyi

    \scriptscriptfont\cyifam=\fivecyi

\def\id{{\mbox{1 \hskip -7pt 1}}}
\newcommand{\sgn}{{\mathit s  \mathit g\mathit  n}}
 \newcommand{\lon}{\longrightarrow}
 \newcommand{\bu}{\bullet}
 
 \newcommand{\rar}{\rightarrow}

 \newcommand{\CP}{{\mathbb C} {\mathbb P}}
 \newcommand{\End}{\mathsf{End}}
\newcommand{\p}{{\partial}}
\newcommand{\Id}{{\mathrm{Id}}}

\newcommand{\assin}{\mathcal A ss_\infty}

\newcommand{\Q}{{\mathbb Q}}
 \newcommand{\Z}{{\mathbb Z}}
 \newcommand{\bS}{{\mathbb S}}
 
 \newcommand{\C}{{\mathbb C}}
 \newcommand{\R}{{\mathbb R}}
 \newcommand{\N}{{\mathbb N}}
 \newcommand{\K}{{\mathbb K}}

 \newcommand{\ot}{\otimes}

\newcommand{\bW}{{\mathbf W}}

\newcommand{\Def}{\mathsf{Def}}

\newcommand{\BV}{{\mathcal B}{\mathcal V}}

\newcommand{\RH}{\cR \caH ra}

 \newcommand{\Beq}{\begin{equation}}
 \newcommand{\Eeq}{\end{equation}}
 \newcommand{\Beqr}{\begin{eqnarray}}
 \newcommand{\Eeqr}{\end{eqnarray}}
 \newcommand{\Beqrn}{\begin{eqnarray*}}
 \newcommand{\Eeqrn}{\end{eqnarray*}}
 \newcommand{\Ba}{\begin{array}}
 \newcommand{\Ea}{\end{array}}
 \newcommand{\Bi}{\begin{itemize}}
 \newcommand{\Ei}{\end{itemize}}
 \newcommand{\Bc}{\begin{center}}
 \newcommand{\Ec}{\end{center}}

 \newcommand{\fg}{{\mathfrak g}}

 \newcommand{\fo}{{\mathfrak o}}

\newcommand{\fs}{{\mathfrak s}}

 \newcommand{\f}{{\mathcal O}}
 
 \newcommand{\cB}{{\mathcal B}}

 \newcommand{\cE}{{\mathcal E}}
 \newcommand{\cF}{{\mathcal F}}
 \newcommand{\cG}{{\mathcal G}}
 \newcommand{\caH}{{\mathcal H}}

 \newcommand{\cM}{{\mathcal M}}
 
 \newcommand{\cP}{{\mathcal P}}
 
 \newcommand{\cR}{{\mathcal R}}

 \newcommand{\cV}{{\mathcal V}}
 \newcommand{\cW}{{\mathcal W}}

 \newcommand{\al}{\alpha}
 \newcommand{\be}{\beta}
 \newcommand{\ga}{\gamma}
 
 \newcommand{\Ga}{\Gamma}

 \newcommand{\la}{\lambda}

 \newcommand{\Img}{{\mathsf I\mathsf m}\, }
 \newcommand{\Hom}{{\mathrm H\mathrm o\mathrm m}}
 
 \newcommand{\sip}{\smallskip}
 \newcommand{\bip}{\bigskip}
 \newcommand{\mip}{\vspace{2.5mm}}

\newcommand{\Gal}{{\mathsf{Gal}}}

\newcommand{\HocBVd}{\mathcal{H}\mathit{o}\mathcal{BV}^{com}_{d}}

\newcommand{\LoB}{\mathcal{L}\mathit{ieb}^\diamond}

\newcommand{\LBd}{\mathcal{L}\mathit{ieb}_{d}}

\newcommand{\HoLBd}{\mathcal{H}\mathit{olieb}_{d}}

\newcommand{\LoBd}{\mathcal{L}\mathit{ieb}_{d}^\diamond}

\newcommand{\HoLoBd}{\mathcal{H}\mathit{olieb}_{d}^\diamond}

\newcommand{\HoLoB}{\mathcal{H}\mathit{olieb}^\diamond}

\theoremstyle{plain}
\swapnumbers

\newtheorem{prop-def}[theorem]{Proposition-definition}

\newtheorem{f-theorem}{Formality Theorem}[section]
\newtheorem{main-theorem}{Main~Theorem}[section]
\newtheorem{section-theorem}{Theorem}[section]

\theoremstyle{definition}

\begin{document}

 \sloppy

 \newenvironment{proo}{\begin{trivlist} \item{\sc {Proof.}}}
  {\hfill $\square$ \end{trivlist}}

\long\def\symbolfootnote[#1]#2{\begingroup%
\def\thefootnote{\fnsymbol{footnote}}\footnote[#1]{#2}\endgroup}

\title{Prop of ribbon hypergraphs and\\ strongly homotopy involutive Lie bialgebras}

\author{Sergei~Merkulov}
\address{Sergei~Merkulov:  Mathematics Research Unit, Luxembourg University,  Grand Duchy of Luxembourg }
\email{sergei.merkulov@uni.lu}

\begin{abstract} For any integer $d$  we introduce  a prop  $\cR\caH ra_d$ of oriented ribbon hypergraphs (in which ``edges" can connect more than two vertices) and prove that it admits a canonical morphism of props,
$$
\HoLoBd \lon \cR\caH ra_d,
$$
 $\HoLoBd$ being   the (degree shifted) minimal resolution of prop of involutive Lie bialgebras,
which is non-trivial on every generator of  $\HoLoBd$. We obtain two applications of this general construction.

\sip

First
we show that for any graded vector space $W$ equipped with a family of cyclically (skew)symmetric higher products,
$$
\Theta_n: (\ot^n W[d])_{\Z_n} \lon \K[1+d], \ \ \ \ \ n\geq 1,
$$
the associated vector space of cyclic words
 $Cyc(W)=\oplus_{n\geq 0} (\ot^n W)_{\Z_n} $ has a combinatorial $\HoLoBd$-structure. As an illustration we construct for each natural number
 $N\geq 1$ an explicit combinatorial strongly homotopy involutive Lie bialgebra structure on the vector space of cyclic words in $N$ graded letters which extends the well-known Schedler's necklace  Lie bialgebra structure from the formality theory of the Goldman-Turaev Lie bialgebra in genus zero.

 \sip

Second, we introduced new (in general, non-trivial) operations in string topology. Given any closed connected and simply connected manifold $M$ of dimension $\geq 4$. We show that the reduced equivariant homology $\bar{H}_\bu^{S^1}(LM)$ of the space $LM$ of free loops in $M$ carries a canonical representation of the dg prop  $\HoLoB_{2-n}$  on  $\bar{H}_\bu^{S^1}(LM)$ controlled by four ribbon hypergraphs explicitly shown in this paper.

\bip


\end{abstract}
 \maketitle
 \markboth{}{}

{\Large
\section{\bf Introduction}
}

\subsection{Involutive Lie bialgebras} Lie bialgebras have been introduced by V.\ Drinfeld in \cite{D1}  in his studies of Yang-Baxter equations and the deformation theory of universal enveloping algebras. Nowadays involutive Lie bialgebras are used in many different areas of mathematics --- in algebra, geometry, string topology, contact topology, theory of moduli spaces of algebraic curves, etc. (see, e.g., articles and books \cite{AKKN, AKKN2, Ba1, CEG, CFL, Ch, ChSu, D1,D2, DCTT, ES, Ma, MW, NW, Tu, Sch} as well as references cited there). The problem  to prove Koszulness of the prop of involutive Lie bialgebras (and hence to understand their their homotopy theory) turned out to be a rather non-trivial one --- it stayed open for almost a decade until its solution in \cite{CMW}. It is a remarkable fact that the deformation theory of that prop is controlled by the mysterious Grothendieck-Teichm\"uller group \cite{MW2} which appears in many different areas of mathematics and explains, perhaps, the diversity of important mathematical problems which involve involutive Lie bialgebras.

\sip

The main result of this paper is an explicit construction of several families of strongly homotopy involutive Lie
bialgebras using a new prop of ribbon hypergraphs $\caH \cG ra_d$, $\forall d\in \Z$. We show two applications of our results  --- one is a purely combinatorial one generalizing a nice construction of T.\ Schedler in \cite{Sch}, the second one deals with  the beautiful theory of string topology introduced by M.\ Chas and D.\ Sullivan in \cite{CS} and introduces four new strongly homotopy operations on the reduced equivariant homology $\bar{H}_\bu^{S^1}(LM)$ of the space of free loops in an arbitrary closed simply connected manifold $M$ of dimension $\geq 4$.

\sip

 This work is much motivated by a particular class of  involutive Lie bialgebra structures
  originating in the theory of punctured Riemann surfaces \cite{G,Tu,AKKN,AKKN2} and in the string topology of an arbitrary closed simply connected manifold  \cite{CS,CEG,CFL,NW}.

  \sip

 Let $\widehat{\K}\langle\pi_1(\Sigma_{0,N+1})\rangle$ stand for the completed group algebra of the fundamental group $\pi_1(\Sigma_{0,N+1}, \K)$ of
  the genus zero Riemann surface $\Sigma_{0,N+1}$ with $N+1$ boundary components, $N\geq 2$, and $H_1(\Sigma_{0,N+1})$ for its first homology group over $\K$. Let
 $$
 \widehat{\fg}[\Sigma_{0,N+1}]:= \frac{\widehat{\K}\langle\pi_1(\Sigma_{0,N+1})\rangle}{[ \widehat{\K}\langle\pi_1(\Sigma_{0,N+1}), \widehat{\K}\langle\pi_1(\Sigma_{0,N+1})\rangle]}\
 $$
be the (completed) vector space spanned over a field $\K$ of characteristic zero by free homotopy classes of loops in $\Sigma_{0,N+1}$. Using intersections and self-intersection of loops Goldman and Turaev \cite{G,Tu}
made this vector space into a filtered involutive Lie bialgebra\footnote{Strictly
 speaking this structure was originally defined modulo constant loops, i.e.\
 on the quotient space  $\widehat{\fg}[\Sigma_{0,N+1}]/\K$, but a choice of
 framing  on  $\Sigma_{0,N+1}$ permits us to extend that structure to
 the whole space. Put another way, the induced Lie bialgebra structure on
 $\widehat{\fg}[\Sigma_{0,N+1}]$ is canonical only on the quotient space
  $Cyc^\bu(W_N)/\K\id$, $\id$ being the empty cyclic word,
 and its extension to the whole space depends on some additional choices.
 }. Let
$$
\mathrm{gr}\widehat{\fg}[\Sigma_{0,N+1}] := \frac{\widehat{\ot^\bu} H_1(\Sigma_{0,N+1},\K)}{[ \widehat{\ot^\bu} H_1(\Sigma_{0,N+1},\K), \widehat{\ot^\bu} H_1(\Sigma_{0,N+1},\K)]}\simeq Cyc(W_N):= \prod_{n\geq 0} (\ot^n W_N)_{\Z_n}
$$
be the {\em associated graded}\, involutive Lie bialgebra where
 $$
 W_N=\text{span}_\K\left\langle x_1,\ldots x_N\right\rangle
 $$
stands for the vector spaces generated by $N$ formal letters $x_1, \ldots, x_N$ (corresponding to the standard generators of $H_1(\Sigma_{0,N+1},\K)$).
The {\em formality theorem}  \cite{AKKN,AN, Ma} establishes a highly non-trivial isomorphism of Lie bialgebras
$$
\widehat{\fg}[\Sigma_{0,N+1}]  \lon \
\mathrm{gr}\widehat{\fg}[\Sigma_{0,N+1}]
$$
 which depends on the choice of a Drinfeld associator. Thus the Goldman-Turaev Lie bialgebra structure can be understood in terms of  its much simpler graded associated version which admits a purely combinatorial description. In fact, it admits two purely combinatorial
 descriptions. The first one is due to the general construction by Schedler \cite{Sch}
 which associates to any quiver a so called necklace Lie bialgebra; the particular involutive Lie bialgebra structure on $\mathrm{gr}\widehat{\fg}[\Sigma_{0,N+1}]\simeq Cyc(W_N)$ is the necklace one corresponding  to the following quiver
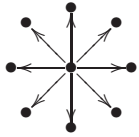
\Beq\label{1: quiver for W_N}
 \xy
 (0,0)*{\bu}="a",
(8,0)*{\bu}="b1",
(6,-6)*{\bu}="b2",
(0,-8)*{\bu}="b3",
(-6,-6)*{\bu}="b4",
(-8,0)*{\bu}="b5",
(-6,6)*{\bu}="b6",
(0,8)*{\bu}="b7",
(6,6)*{\bu}="b8",
\ar @{->} "a";"b1" <0pt>
\ar @{->} "a";"b2" <0pt>
\ar @{->} "a";"b3" <0pt>
\ar @{->} "a";"b4" <0pt>
\ar @{->} "a";"b5" <0pt>
\ar @{->} "a";"b6" <0pt>
\ar @{->} "a";"b7" <0pt>
\ar @{->} "a";"b8" <0pt>
\endxy
\Eeq
with $N$ legs.
Put another way, for any natural number $N\geq 2$ T.\ Schedler's construction\footnote{It is worth emphasizing that Schedler's construction depends on the choice of a basis  $(x_1,\ldots, x_N)$ in $W_N$, i.e.\ it depends essentially only on the natural number $N$.} gives us an involutive Lie bialgebra structure on $Cyc(W_N)$ which admits a nice geometric interpretation.

\sip

The second combinatorial description of the necklace Lie bialgebra structure on $Cyc(W_N)$ involves the $d=1$ case of a family of props of ribbon graphs $\cR \cG ra_d$, $d\in \Z$, which come equipped with canonical morphisms \cite{MW}
$$
\rho: \LoBd \lon \cR \cG ra_d
$$
from the prop of involutive Lie  bialgebras; the map $\rho$ is non-trivial on both Lie and coLie generators of $\LoBd$  (which are assigned homological degree $1-d$ so that the case $d=1$ corresponds to the ordinary Lie bialgebras). It was shown in \cite{MW} that for any $\Z$-graded vector space $V$ equipped with a pairing (which for $d=1$ is nothing but a skew-symmetric scalar product in $W$)
$$
\Ba{rccc}
\Theta_2: & \odot^2(W[d]) & \lon & \K[1+d]
\Ea
$$
there is an associated representation
$$
\rho_{\Theta_2}: \cR \cG ra_d \lon \cE nd_{Cyc(W)}
$$
of the prop of ribbon graphs in the vector space
$$
 Cyc(W):= \prod_{n\geq 1} (\ot^n W)_{\Z_n}
$$
spanned by cyclic words in elements of $W$, and hence there is an induced via the composition $\rho_{\Theta_2}\circ \rho$ an involutive  Lie bialgebra structure in $Cyc(W)$. The vector space $W_N$ has no natural pairings so we can not apply this construction immediately to get the necklace Lie bialgebra structure on $W_N$. However a certain ``doubling" trick explained in \S {\ref{2: Sch liab}} does the job and gives us  a canonical representation
\Beq\label{1:  rho_N}
\rho_N: \cR \cG ra _1 \lon \cE nd_{Cyc(W_N)}
\Eeq
which induces via the composition $\rho_N \circ \rho$ the required necklace Lie bialgebra structure on $Cyc(W_N)$ (modulo terms depending on a particular choice
of framing on $\Sigma_{0,N+1}$, i.e.\ both structures fully agree on the quotient
space $Cyc(W_N)/\K\id$, where $\id$ stands for the empty cyclic word; in fact, the particular map $\rho_N$ we construct in  \S {\ref{2: Sch liab}} corresponds to the blackboard framing as it was explained to the author by Yusuke Kuno \cite{Ku}.

\sip

Given any Poincare duality algebra $A=\K \oplus \bar{A}$ in degree $n$  (say, a one which models a closed simply connected $n$-dimensional manifold $M$), the associated space $W= \bar{A}^*[-1]$ comes equipped with a non-degenerate pairing $\Theta_2$ as above with $d=3-n$. Hence the above construction gives us a representation of $\cR \cG ra_{3-n}$ in $Cyc(W)$ and hence induces a $\LoB_{3-n}$-algebra structure on that space. It has been proven in \cite{CEG} that $\LoB_{3-n}$-operations on $Cyc(\bar{A}^*[-1])$ respect the Hochschild differential $d_H$ and hence induce a  $\LoB_{3-n}$-algebra structure on its
cohomology $H^\bu(Cyc(\bar{A}^*[-1]), d_H)$. If $A$ happens to be a model of a closed simply connected $n$ manifold $M$, then there is an isomorphism of (co)homology group
$$
\bar{H}_\bu^{S^1}(LM) \stackrel{\simeq}{\lon} H^\bu(Cyc(\bar{A}^*[-1]), d_H)
$$
where the l.h.s.\ stands for the reduced (negatively graded) equivariant cohomology of the space $LM$ of free loops in $M$.
Hence $\bar{H}_\bu^{S^1}(LM)$ comes equipped with a $\LoB_{3-n}$-algebra structure which was first discovered in a purely geometric way by M.\ Chas and D.\ Sullivan  in \cite{CS}.

\bip

In this paper we introduce, for each integer $d\in \Z$, a  {\em prop of ribbon hypergraphs $\caH \cG ra_d$} and show that there is a morphism
of dg props
$$
\rho^\diamond: \HoLoBd \lon \caH \cG ra_d
$$
which is {\em non-trivial on every generator of}\, the minimal resolution $\HoLoBd$  of the prop $\LoBd$. There is a natural commutative diagram
$$
\xymatrix{
\HoLoBd\ar[d]_p \ar[r]^{\rho^\diamond}    &   \caH \cG ra_d \ar[d]^q\\
\LoBd  \ar[r]^{\rho} &    \cR \cG ra_d
}
$$
where the left vertical arrow $p$ is the canonical quasi-isomorphism and the right arrow $q$ is a ``forgetful" map sending to zero all ribbon hypergraphs
with at least one non-bivalent hyperedge.

\sip

Given any graded vector space $W$ equipped with cyclically (skew)invariant maps
\Beq\label{1: Theta_n}
\Ba{rccc}
\Theta_n: & \odot^n(W[d])_{\Z_n} & \lon & \K[1+d], \ \ \ \ n\in \N_{\geq 1},
\Ea
\Eeq
there is a canonical representation
$$
\rho_{\Theta_{\bu}}: \caH \cG ra_d \lon \cE nd_{Cyc(W)}
$$
of the prop of ribbon hypergraphs, and hence a canonical strongly homotopy
involutive Lie bialgebra structure
$$
\rho_{\Theta_{\bu}}\circ \rho^\diamond: \HoLoBd \lon \cE nd_{Cyc(W)}
$$
on the space of cyclic words $Cyc(W)$ in letters from $W$. If all $\Theta_n$ vanish except for $n=2$ we recover the
previous result from \cite{MW}.

\sip

We consider two applications of the above construction, one is in purely algebra and another one is in string topology.

\sip

Using multi-tuple generalization of the ``doubling" trick
used in the construction of representation (\ref{1:  rho_N}) we introduce an
explicit highly non-trivial strongly homotopy involutive Lie bialgebra structure,
$$
\hat{\rho}_N: \HoLoB_1 \lon \cE nd_{Cyc(\widehat{W}_N)}
$$
in the vector space generated by cyclic words in $\Z$-graded letters
$\{x_1[-p], \ldots, x_N[-p]\}_{p\in \N}$, where $x_i[-p]$ stands for the formal letter $x_i$ to which
we assigned homological degree $p$. When all letters are concentrated in degree zero, one recovers
Schedler's necklace Lie bialgebra associated to the quiver (\ref{1: quiver for W_N}).

\sip

Let $M$ be an arbitrary closed simply connected $n$-dimensional manifold, and $\Omega_M^\bu$ its de Rham algebra. Using the integration morphism (more precisely, its version in a Poincare duality algebra modelling $M$)
$$
\Ba{ccc}
\ot^3 \Omega_M^\bu & \lon & \R\\
  \al\ot \be\ot\ga & \lon & \int_M \al\wedge \be\wedge \ga
\Ea
$$
we prove that the reduced equivariant cohomology $\bar{H}_\bu^{S^1}(LM)$
  of the space $LM$ of free loops in $M$
comes equipped canonically with a  $\HoLoB_{2-n}$-structure
whose  only possibly non-trivial operations are controlled by the following four generators of $\HoLoB_{2-n}$ (see \S 2 for their definition),
$$
\Ba{c}\resizebox{8mm}{!}{\xy
(-5,-6)*{};
(0,0)*+{_0}*\cir{}
**\dir{-};
(0,-6)*{};
(0,0)*+{_0}*\cir{}
**\dir{-};
(5,-6)*{};
(0,0)*+{_0}*\cir{}
**\dir{-};
(0,6)*{};
(0,0)*+{_0}*\cir{}
**\dir{-};
(0,8)*{_1};
(-5,-8)*{_1};
(0,-8)*{_2};
(5,-8)*{_3};
\endxy}\Ea\ \ , \ \
\Ba{c}\resizebox{8mm}{!}{\xy
(-5,6)*{};
(0,0)*+{_0}*\cir{}
**\dir{-};
(0,6)*{};
(0,0)*+{_0}*\cir{}
**\dir{-};
(5,6)*{};
(0,0)*+{_0}*\cir{}
**\dir{-};
(0,-6)*{};
(0,0)*+{_0}*\cir{}
**\dir{-};
(0,-8)*{_1};
(-5,8)*{_1};
(0,8)*{_2};
(5,8)*{_3};
\endxy}\Ea\ \ , \ \
\Ba{c}\resizebox{7mm}{!}{\xy
(-5,-6)*{};
(0,0)*+{_0}*\cir{}
**\dir{-};
(5,6)*{};
(0,0)*+{_0}*\cir{}
**\dir{-};
(5,-6)*{};
(0,0)*+{_0}*\cir{}
**\dir{-};
(-5,6)*{};
(0,0)*+{_0}*\cir{}
**\dir{-};
(-5,8)*{_1};
(5,8)*{_2};
(-5,-8)*{_1};
(5,-8)*{_2};
\endxy}\Ea\ \ , \ \
\Ba{c}\resizebox{2.5mm}{!}{\xy
(0,-6)*{};
(0,0)*+{_1}*\cir{}
**\dir{-};
(0,6)*{};
(0,0)*+{_1}*\cir{}
**\dir{-};
(0,8)*{_1};
(0,-8)*{_1};
\endxy}\Ea
$$
which in turn are respectively controlled via the above mentioned morphism $\rho^\diamond$ by the following ribbon hypergraphs from  $\RH_{2-n}$,
$$
\xy
(0,4)*{\ast}="1";
(-5,-4)*{\circ}="2";
(0,-4)*{\circ}="3";
(5,-4)*{\circ}="4";
\ar @{-} "1";"2" <0pt>
\ar @{-} "1";"3" <0pt>
\ar @{-} "1";"4" <0pt>
\endxy\ \ ,  \
\ \ \ \xy
(0,5)*{\ast}="1";
(0,-4)*{\circ}="3";
"1";"3" **\crv{(4,0) & (4,1)};
"1";"3" **\crv{(-4,0) & (-4,-1)};
\ar @{-} "1";"3" <0pt>
\endxy\ \ ,  \ \
\xy
(0,5)*{\ast}="1";
(-5,-4)*{\circ}="3";
(5,-4)*{\circ}="4";
"1";"3" **\crv{(-3,5) & (5,4)};
"1";"3" **\crv{(-5,2) & (-5,2)};
\ar @{-} "1";"4" <0pt>
\endxy\ \ ,
\ \ \ \xy
(0,5)*{\ast}="1";
(0,-4)*{\circ}="3";
"1";"3" **\crv{(-5,2) & (5,2)};
"1";"3" **\crv{(5,2) & (-5,2)};
"1";"3" **\crv{(-7,7) & (-7,-7)};
\endxy
%
%
%
%
$$
with vertices and boundaries appropriately (skew)symmetrized (so that their labelling is omitted).
It is easy to see that these four new $\HoLoB_{2-n}$ operations on $\bar{H}_\bu^{S^1}(LM)$ are non-trivial for $M=\CP^n$ with $n\geq 3$.

\subsection{Structure of the paper}
In \S 2 we remind basic facts about the prop of (degree shifted) involutive Lie bialgebras $\LoBd$ and its minimal resolution $\HoLoBd$, and their interrelations with
$\caH o\cB\cV_d^{com}$-algebras (see \cite{CMW} for  more details and proofs). In \S 3 we describe in detail the prop of ribbon hypergraphs which is a natural extension of the prop of ribbon graphs introduced in \cite{MW}, and then construct its canonical representation in the space of cyclic words $Cyc(W)$ associated with any graded vector space $W$ equipped with a family of cyclically (skew)symmetric higher products (which is a rather straightforward generalization of the similar construction in \cite{MW}). In \S 4 we explain the main  (and not that straightforward) result of this paper --- a
 construction of an explicit morphism of props,
$\HoLoBd\rar \RH_d$, which is {\em non-trivial on every generator of $\HoLoBd$}. This result  gives us a large family of explicit strongly homotopy involutive Lie bialgebras; in particular, for any natural number $N\in \N_{\geq 1}$ we show in \S 5  an explicit strongly homotopy involutive Lie bialgebra structure on the vector space of cyclic words in $\Z$-graded formal letters which extends the well-known Schedler's necklace  Lie bialgebra structure \cite{Sch}  from the formality theory of the Goldman-Turaev Lie bialgebra in genus zero. In \S 6 we discuss an application of the main construction to string topology.

\subsection{Some notation}
 The set $\{1,2, \ldots, n\}$ is abbreviated to $[n]$;  its group of automorphisms is
denoted by $\bS_n$. The trivial (resp., sign) one-dimensional representation of
 $\bS_n$ is denoted by $\id_n$ (resp.,  $\sgn_n$).
 The cardinality of a finite set $A$ is denoted by $\# A$.
If $V=\oplus_{i\in \Z} V^i$ is a graded vector space, then
$V[n]$ stands for the graded vector space with $V[n]^i:=V^{i+n}$; for $v\in V^i$ we set $|v|:=i$. The canonical degree $-1$ isomorphism $V\rar V[1]$ is denoted by $\fs$;
one has $|\fs^k v|=|v| -k$ for any homogeneous element $v\in V$.

\sip

For a
prop(erad) $\cP$ we denote by $\cP\{n\}$ a prop(erad) which is uniquely defined by
 the following property:
for any graded vector space $W$ a representation
of $\cP\{n\}$ in $W$ is identical to a representation of  $\cP$ in $W[n]$.

\sip

For a  module $V$ over a group $G$ we denote by $V_G$
 the vector space of coinvariants:
$V/\{g(v) - v\ |\ v\in V, g\in G\}$ and by $V^G$  the subspace
of invariants: $\{\forall g\in  G\ :\  g(v)=v,\ v\in V\}$. We always work over a field $\K$ of characteristic zero so that if $G$ is finite, then these
spaces are canonically isomorphic, $V_G \cong V^G$.



\mip

{\bf Acknowledgement}. It is a pleasure to thank Anton Alekseev, Martin Kassabov, Nariya Kawazumi, Anton Khoroshkin, Yusuke Kuno,  Florian Naef and Thomas Willwacher  for valuable discussions. This text is based partially on the talk given by the author
 at the workshop ``Poisson geometry of moduli spaces, associators and quantum field theory" organized
 at the Simons Center of Geometry and Physics in June 2018. Most of this text was also written there, and the author acknowledges with thanks the excellent working conditions at the SCGP.

\bip

\bip

{\large
\section{\bf  Strongly homotopy involutive Lie bialgebras  as \\
commutative $BV$ algebras and vice versa}
}

\bip

\subsection{\bf Reminder on the prop of strongly homotopy (involutive) Lie bialgebras \cite{CMW}}\label{2: Reminder on Liebs} Let
$$
\LBd:=\cF ree\langle E\rangle/\langle\cR\rangle,
$$
be the quotient
of the free prop(erad) generated by an  $\bS$-bimodule $E=\{E(m,n)\}_{m,n\geq 1}$ with
$$
E(m,n):=\left\{
\Ba{ll} \id_1\ot \sgn_2^{\ot d}[d-1]=\mbox{span}\left\langle
\Ba{c}\begin{xy}
 <0mm,-0.55mm>*{};<0mm,-2.5mm>*{}**@{-},
 <0.5mm,0.5mm>*{};<2.2mm,2.2mm>*{}**@{-},
 <-0.48mm,0.48mm>*{};<-2.2mm,2.2mm>*{}**@{-},
 <0mm,0mm>*{\circ};<0mm,0mm>*{}**@{},
 <0.5mm,0.5mm>*{};<2.7mm,2.8mm>*{^{_2}}**@{},
 <-0.48mm,0.48mm>*{};<-2.7mm,2.8mm>*{^{_1}}**@{},
 \end{xy}\Ea
=(-1)^{d}
\Ba{c}\begin{xy}
 <0mm,-0.55mm>*{};<0mm,-2.5mm>*{}**@{-},
 <0.5mm,0.5mm>*{};<2.2mm,2.2mm>*{}**@{-},
 <-0.48mm,0.48mm>*{};<-2.2mm,2.2mm>*{}**@{-},
 <0mm,0mm>*{\circ};<0mm,0mm>*{}**@{},
 <0.5mm,0.5mm>*{};<2.7mm,2.8mm>*{^{_1}}**@{},
 <-0.48mm,0.48mm>*{};<-2.7mm,2.8mm>*{^{_2}}**@{},
 \end{xy}\Ea
   \right\rangle & \text{if}\ m=2,n=1
 \\
 \sgn_2^{\ot d}\ot \id_1[d-1]=\mbox{span}\left\langle
\Ba{c}\begin{xy}
 <0mm,0.66mm>*{};<0mm,3mm>*{}**@{-},
 <0.39mm,-0.39mm>*{};<2.2mm,-2.2mm>*{}**@{-},
 <-0.35mm,-0.35mm>*{};<-2.2mm,-2.2mm>*{}**@{-},
 <0mm,0mm>*{\circ};<0mm,0mm>*{}**@{},
   <0.39mm,-0.39mm>*{};<2.9mm,-4mm>*{^{_2}}**@{},
   <-0.35mm,-0.35mm>*{};<-2.8mm,-4mm>*{^{_1}}**@{},
\end{xy}\Ea
=(-1)^{d}
\Ba{c}\begin{xy}
 <0mm,0.66mm>*{};<0mm,3mm>*{}**@{-},
 <0.39mm,-0.39mm>*{};<2.2mm,-2.2mm>*{}**@{-},
 <-0.35mm,-0.35mm>*{};<-2.2mm,-2.2mm>*{}**@{-},
 <0mm,0mm>*{\circ};<0mm,0mm>*{}**@{},
   <0.39mm,-0.39mm>*{};<2.9mm,-4mm>*{^{_1}}**@{},
   <-0.35mm,-0.35mm>*{};<-2.8mm,-4mm>*{^{_2}}**@{},
\end{xy}\Ea
\right\rangle
& \text{if}\ m=1,n=2\\
0 & \text{otherwise}
\Ea
\right.
$$
by the ideal $\langle\cR\rangle$ generated by the following relations
\Beq\label{2: R for LieB}
\cR:\left\{
\Ba{c}
\Ba{c}\resizebox{7mm}{!}{
\begin{xy}
 <0mm,0mm>*{\circ};<0mm,0mm>*{}**@{},
 <0mm,-0.49mm>*{};<0mm,-3.0mm>*{}**@{-},
 <0.49mm,0.49mm>*{};<1.9mm,1.9mm>*{}**@{-},
 <-0.5mm,0.5mm>*{};<-1.9mm,1.9mm>*{}**@{-},
 <-2.3mm,2.3mm>*{\circ};<-2.3mm,2.3mm>*{}**@{},
 <-1.8mm,2.8mm>*{};<0mm,4.9mm>*{}**@{-},
 <-2.8mm,2.9mm>*{};<-4.6mm,4.9mm>*{}**@{-},
   <0.49mm,0.49mm>*{};<2.7mm,2.3mm>*{^3}**@{},
   <-1.8mm,2.8mm>*{};<0.4mm,5.3mm>*{^2}**@{},
   <-2.8mm,2.9mm>*{};<-5.1mm,5.3mm>*{^1}**@{},
 \end{xy}}\Ea
 +
\Ba{c}\resizebox{7mm}{!}{\begin{xy}
 <0mm,0mm>*{\circ};<0mm,0mm>*{}**@{},
 <0mm,-0.49mm>*{};<0mm,-3.0mm>*{}**@{-},
 <0.49mm,0.49mm>*{};<1.9mm,1.9mm>*{}**@{-},
 <-0.5mm,0.5mm>*{};<-1.9mm,1.9mm>*{}**@{-},
 <-2.3mm,2.3mm>*{\circ};<-2.3mm,2.3mm>*{}**@{},
 <-1.8mm,2.8mm>*{};<0mm,4.9mm>*{}**@{-},
 <-2.8mm,2.9mm>*{};<-4.6mm,4.9mm>*{}**@{-},
   <0.49mm,0.49mm>*{};<2.7mm,2.3mm>*{^2}**@{},
   <-1.8mm,2.8mm>*{};<0.4mm,5.3mm>*{^1}**@{},
   <-2.8mm,2.9mm>*{};<-5.1mm,5.3mm>*{^3}**@{},
 \end{xy}}\Ea
 +
\Ba{c}\resizebox{7mm}{!}{\begin{xy}
 <0mm,0mm>*{\circ};<0mm,0mm>*{}**@{},
 <0mm,-0.49mm>*{};<0mm,-3.0mm>*{}**@{-},
 <0.49mm,0.49mm>*{};<1.9mm,1.9mm>*{}**@{-},
 <-0.5mm,0.5mm>*{};<-1.9mm,1.9mm>*{}**@{-},
 <-2.3mm,2.3mm>*{\circ};<-2.3mm,2.3mm>*{}**@{},
 <-1.8mm,2.8mm>*{};<0mm,4.9mm>*{}**@{-},
 <-2.8mm,2.9mm>*{};<-4.6mm,4.9mm>*{}**@{-},
   <0.49mm,0.49mm>*{};<2.7mm,2.3mm>*{^1}**@{},
   <-1.8mm,2.8mm>*{};<0.4mm,5.3mm>*{^3}**@{},
   <-2.8mm,2.9mm>*{};<-5.1mm,5.3mm>*{^2}**@{},
 \end{xy}}\Ea
 \ \ , \ \
\Ba{c}\resizebox{8.4mm}{!}{ \begin{xy}
 <0mm,0mm>*{\circ};<0mm,0mm>*{}**@{},
 <0mm,0.69mm>*{};<0mm,3.0mm>*{}**@{-},
 <0.39mm,-0.39mm>*{};<2.4mm,-2.4mm>*{}**@{-},
 <-0.35mm,-0.35mm>*{};<-1.9mm,-1.9mm>*{}**@{-},
 <-2.4mm,-2.4mm>*{\circ};<-2.4mm,-2.4mm>*{}**@{},
 <-2.0mm,-2.8mm>*{};<0mm,-4.9mm>*{}**@{-},
 <-2.8mm,-2.9mm>*{};<-4.7mm,-4.9mm>*{}**@{-},
    <0.39mm,-0.39mm>*{};<3.3mm,-4.0mm>*{^3}**@{},
    <-2.0mm,-2.8mm>*{};<0.5mm,-6.7mm>*{^2}**@{},
    <-2.8mm,-2.9mm>*{};<-5.2mm,-6.7mm>*{^1}**@{},
 \end{xy}}\Ea
 +
\Ba{c}\resizebox{8.4mm}{!}{ \begin{xy}
 <0mm,0mm>*{\circ};<0mm,0mm>*{}**@{},
 <0mm,0.69mm>*{};<0mm,3.0mm>*{}**@{-},
 <0.39mm,-0.39mm>*{};<2.4mm,-2.4mm>*{}**@{-},
 <-0.35mm,-0.35mm>*{};<-1.9mm,-1.9mm>*{}**@{-},
 <-2.4mm,-2.4mm>*{\circ};<-2.4mm,-2.4mm>*{}**@{},
 <-2.0mm,-2.8mm>*{};<0mm,-4.9mm>*{}**@{-},
 <-2.8mm,-2.9mm>*{};<-4.7mm,-4.9mm>*{}**@{-},
    <0.39mm,-0.39mm>*{};<3.3mm,-4.0mm>*{^2}**@{},
    <-2.0mm,-2.8mm>*{};<0.5mm,-6.7mm>*{^1}**@{},
    <-2.8mm,-2.9mm>*{};<-5.2mm,-6.7mm>*{^3}**@{},
 \end{xy}}\Ea
 +
\Ba{c}\resizebox{8.4mm}{!}{ \begin{xy}
 <0mm,0mm>*{\circ};<0mm,0mm>*{}**@{},
 <0mm,0.69mm>*{};<0mm,3.0mm>*{}**@{-},
 <0.39mm,-0.39mm>*{};<2.4mm,-2.4mm>*{}**@{-},
 <-0.35mm,-0.35mm>*{};<-1.9mm,-1.9mm>*{}**@{-},
 <-2.4mm,-2.4mm>*{\circ};<-2.4mm,-2.4mm>*{}**@{},
 <-2.0mm,-2.8mm>*{};<0mm,-4.9mm>*{}**@{-},
 <-2.8mm,-2.9mm>*{};<-4.7mm,-4.9mm>*{}**@{-},
    <0.39mm,-0.39mm>*{};<3.3mm,-4.0mm>*{^1}**@{},
    <-2.0mm,-2.8mm>*{};<0.5mm,-6.7mm>*{^3}**@{},
    <-2.8mm,-2.9mm>*{};<-5.2mm,-6.7mm>*{^2}**@{},
 \end{xy}}\Ea
 \\
(-1)^d \Ba{c}\resizebox{5mm}{!}{\begin{xy}
 <0mm,2.47mm>*{};<0mm,0.12mm>*{}**@{-},
 <0.5mm,3.5mm>*{};<2.2mm,5.2mm>*{}**@{-},
 <-0.48mm,3.48mm>*{};<-2.2mm,5.2mm>*{}**@{-},
 <0mm,3mm>*{\circ};<0mm,3mm>*{}**@{},
  <0mm,-0.8mm>*{\circ};<0mm,-0.8mm>*{}**@{},
<-0.39mm,-1.2mm>*{};<-2.2mm,-3.5mm>*{}**@{-},
 <0.39mm,-1.2mm>*{};<2.2mm,-3.5mm>*{}**@{-},
     <0.5mm,3.5mm>*{};<2.8mm,5.7mm>*{^2}**@{},
     <-0.48mm,3.48mm>*{};<-2.8mm,5.7mm>*{^1}**@{},
   <0mm,-0.8mm>*{};<-2.7mm,-5.2mm>*{^1}**@{},
   <0mm,-0.8mm>*{};<2.7mm,-5.2mm>*{^2}**@{},
\end{xy}}\Ea
  -
\Ba{c}\resizebox{7mm}{!}{\begin{xy}
 <0mm,-1.3mm>*{};<0mm,-3.5mm>*{}**@{-},
 <0.38mm,-0.2mm>*{};<2.0mm,2.0mm>*{}**@{-},
 <-0.38mm,-0.2mm>*{};<-2.2mm,2.2mm>*{}**@{-},
<0mm,-0.8mm>*{\circ};<0mm,0.8mm>*{}**@{},
 <2.4mm,2.4mm>*{\circ};<2.4mm,2.4mm>*{}**@{},
 <2.77mm,2.0mm>*{};<4.4mm,-0.8mm>*{}**@{-},
 <2.4mm,3mm>*{};<2.4mm,5.2mm>*{}**@{-},
     <0mm,-1.3mm>*{};<0mm,-5.3mm>*{^1}**@{},
     <2.5mm,2.3mm>*{};<5.1mm,-2.6mm>*{^2}**@{},
    <2.4mm,2.5mm>*{};<2.4mm,5.7mm>*{^2}**@{},
    <-0.38mm,-0.2mm>*{};<-2.8mm,2.5mm>*{^1}**@{},
    \end{xy}}\Ea
  - (-1)^{d}
\Ba{c}\resizebox{7mm}{!}{\begin{xy}
 <0mm,-1.3mm>*{};<0mm,-3.5mm>*{}**@{-},
 <0.38mm,-0.2mm>*{};<2.0mm,2.0mm>*{}**@{-},
 <-0.38mm,-0.2mm>*{};<-2.2mm,2.2mm>*{}**@{-},
<0mm,-0.8mm>*{\circ};<0mm,0.8mm>*{}**@{},
 <2.4mm,2.4mm>*{\circ};<2.4mm,2.4mm>*{}**@{},
 <2.77mm,2.0mm>*{};<4.4mm,-0.8mm>*{}**@{-},
 <2.4mm,3mm>*{};<2.4mm,5.2mm>*{}**@{-},
     <0mm,-1.3mm>*{};<0mm,-5.3mm>*{^2}**@{},
     <2.5mm,2.3mm>*{};<5.1mm,-2.6mm>*{^1}**@{},
    <2.4mm,2.5mm>*{};<2.4mm,5.7mm>*{^2}**@{},
    <-0.38mm,-0.2mm>*{};<-2.8mm,2.5mm>*{^1}**@{},
    \end{xy}}\Ea
  -
\Ba{c}\resizebox{7mm}{!}{\begin{xy}
 <0mm,-1.3mm>*{};<0mm,-3.5mm>*{}**@{-},
 <0.38mm,-0.2mm>*{};<2.0mm,2.0mm>*{}**@{-},
 <-0.38mm,-0.2mm>*{};<-2.2mm,2.2mm>*{}**@{-},
<0mm,-0.8mm>*{\circ};<0mm,0.8mm>*{}**@{},
 <2.4mm,2.4mm>*{\circ};<2.4mm,2.4mm>*{}**@{},
 <2.77mm,2.0mm>*{};<4.4mm,-0.8mm>*{}**@{-},
 <2.4mm,3mm>*{};<2.4mm,5.2mm>*{}**@{-},
     <0mm,-1.3mm>*{};<0mm,-5.3mm>*{^2}**@{},
     <2.5mm,2.3mm>*{};<5.1mm,-2.6mm>*{^1}**@{},
    <2.4mm,2.5mm>*{};<2.4mm,5.7mm>*{^1}**@{},
    <-0.38mm,-0.2mm>*{};<-2.8mm,2.5mm>*{^2}**@{},
    \end{xy}}\Ea
 - (-1)^{d}
\Ba{c}\resizebox{7mm}{!}{\begin{xy}
 <0mm,-1.3mm>*{};<0mm,-3.5mm>*{}**@{-},
 <0.38mm,-0.2mm>*{};<2.0mm,2.0mm>*{}**@{-},
 <-0.38mm,-0.2mm>*{};<-2.2mm,2.2mm>*{}**@{-},
<0mm,-0.8mm>*{\circ};<0mm,0.8mm>*{}**@{},
 <2.4mm,2.4mm>*{\circ};<2.4mm,2.4mm>*{}**@{},
 <2.77mm,2.0mm>*{};<4.4mm,-0.8mm>*{}**@{-},
 <2.4mm,3mm>*{};<2.4mm,5.2mm>*{}**@{-},
     <0mm,-1.3mm>*{};<0mm,-5.3mm>*{^1}**@{},
     <2.5mm,2.3mm>*{};<5.1mm,-2.6mm>*{^2}**@{},
    <2.4mm,2.5mm>*{};<2.4mm,5.7mm>*{^1}**@{},
    <-0.38mm,-0.2mm>*{};<-2.8mm,2.5mm>*{^2}**@{},
    \end{xy}}\Ea
    \Ea
\right.
\Eeq
It is called the prop of (degree shifted) {\em Lie bialgebras}.

\sip

The prop of {\em involutive}\, Lie $d$-bialgebras is defined similarly,
$$
\LoBd:=\cF ree\langle E\rangle/\langle\cR_\diamond\rangle
$$
but with a larger set of relations,
$$
\cR_\diamond:= \cR \ \bigsqcup
\Ba{c}\resizebox{4mm}{!}
{\xy
 (0,0)*{\circ}="a",
(0,6)*{\circ}="b",
(3,3)*{}="c",
(-3,3)*{}="d",
 (0,9)*{}="b'",
(0,-3)*{}="a'",
\ar@{-} "a";"c" <0pt>
\ar @{-} "a";"d" <0pt>
\ar @{-} "a";"a'" <0pt>
\ar @{-} "b";"c" <0pt>
\ar @{-} "b";"d" <0pt>
\ar @{-} "b";"b'" <0pt>
\endxy}
\Ea
$$
A representation $\rho: \LBd\rar \cE nd_V$ (resp., $\rho: \LoBd\rar \cE nd_V$) in a graded vector space $V$
provides the latter with two operations
$$
[\ ,\ ]:=\rho\left(\begin{xy}
 <0mm,0.66mm>*{};<0mm,3mm>*{}**@{-},
 <0.39mm,-0.39mm>*{};<2.2mm,-2.2mm>*{}**@{-},
 <-0.35mm,-0.35mm>*{};<-2.2mm,-2.2mm>*{}**@{-},
 <0mm,0mm>*{\circ};<0mm,0mm>*{}**@{},
 \end{xy}\right): \odot^2(V[d])\rar V[1+d], \ \ \ \ \
 \vartriangle:= \left(
 \begin{xy}
 <0mm,-0.55mm>*{};<0mm,-2.5mm>*{}**@{-},
 <0.5mm,0.5mm>*{};<2.2mm,2.2mm>*{}**@{-},
 <-0.48mm,0.48mm>*{};<-2.2mm,2.2mm>*{}**@{-},
 <0mm,0mm>*{\circ};<0mm,0mm>*{}**@{},
 \end{xy}\right): V[d]\lon \odot^2(V[d])[1-2d]
$$
which satisfy the compatibility conditions controlled by the relations
$\cR$ (resp.\ $\cR^\diamond$). If $d=1$, it is precisely the prop of ordinary {\em involutive Lie bialgebras}\, and is often denoted by $\LoB$.

\sip

The properads behind the props $\LBd$ and $\LoBd$ are Koszul so that their minimal resolutions, $\HoLBd$ and respectively $\HoLoBd$, are relatively ``small" (see \cite{CMW,MaVo,Va} and references cited there). The dg prop $\HoLBd$
is generated by the (skew)symmetric corollas of homological degree $1-d(m+n-2)$
$$
\Ba{c}\resizebox{17mm}{!}{\begin{xy}
 <0mm,0mm>*{\circ};<0mm,0mm>*{}**@{},
 <-0.6mm,0.44mm>*{};<-8mm,5mm>*{}**@{-},
 <-0.4mm,0.7mm>*{};<-4.5mm,5mm>*{}**@{-},
 <0mm,0mm>*{};<1mm,5mm>*{\ldots}**@{},
 <0.4mm,0.7mm>*{};<4.5mm,5mm>*{}**@{-},
 <0.6mm,0.44mm>*{};<8mm,5mm>*{}**@{-},
   <0mm,0mm>*{};<-10.5mm,5.9mm>*{^{\sigma(1)}}**@{},
   <0mm,0mm>*{};<-4mm,5.9mm>*{^{\sigma(2)}}**@{},
   <0mm,0mm>*{};<10.0mm,5.9mm>*{^{\sigma(m)}}**@{},
 <-0.6mm,-0.44mm>*{};<-8mm,-5mm>*{}**@{-},
 <-0.4mm,-0.7mm>*{};<-4.5mm,-5mm>*{}**@{-},
 <0mm,0mm>*{};<1mm,-5mm>*{\ldots}**@{},
 <0.4mm,-0.7mm>*{};<4.5mm,-5mm>*{}**@{-},
 <0.6mm,-0.44mm>*{};<8mm,-5mm>*{}**@{-},
   <0mm,0mm>*{};<-10.5mm,-6.9mm>*{^{\tau(1)}}**@{},
   <0mm,0mm>*{};<-4mm,-6.9mm>*{^{\tau(2)}}**@{},
   <0mm,0mm>*{};<10.0mm,-6.9mm>*{^{\tau(n)}}**@{},
 \end{xy}}\Ea
=(-1)^{d(\sigma+\tau)}
\Ba{c}\resizebox{14mm}{!}{\begin{xy}
 <0mm,0mm>*{\circ};<0mm,0mm>*{}**@{},
 <-0.6mm,0.44mm>*{};<-8mm,5mm>*{}**@{-},
 <-0.4mm,0.7mm>*{};<-4.5mm,5mm>*{}**@{-},
 <0mm,0mm>*{};<-1mm,5mm>*{\ldots}**@{},
 <0.4mm,0.7mm>*{};<4.5mm,5mm>*{}**@{-},
 <0.6mm,0.44mm>*{};<8mm,5mm>*{}**@{-},
   <0mm,0mm>*{};<-8.5mm,5.5mm>*{^1}**@{},
   <0mm,0mm>*{};<-5mm,5.5mm>*{^2}**@{},
   <0mm,0mm>*{};<4.5mm,5.5mm>*{^{m\hspace{-0.5mm}-\hspace{-0.5mm}1}}**@{},
   <0mm,0mm>*{};<9.0mm,5.5mm>*{^m}**@{},
 <-0.6mm,-0.44mm>*{};<-8mm,-5mm>*{}**@{-},
 <-0.4mm,-0.7mm>*{};<-4.5mm,-5mm>*{}**@{-},
 <0mm,0mm>*{};<-1mm,-5mm>*{\ldots}**@{},
 <0.4mm,-0.7mm>*{};<4.5mm,-5mm>*{}**@{-},
 <0.6mm,-0.44mm>*{};<8mm,-5mm>*{}**@{-},
   <0mm,0mm>*{};<-8.5mm,-6.9mm>*{^1}**@{},
   <0mm,0mm>*{};<-5mm,-6.9mm>*{^2}**@{},
   <0mm,0mm>*{};<4.5mm,-6.9mm>*{^{n\hspace{-0.5mm}-\hspace{-0.5mm}1}}**@{},
   <0mm,0mm>*{};<9.0mm,-6.9mm>*{^n}**@{},
 \end{xy}}\Ea \ \ \forall \sigma\in \bS_m, \forall\tau\in \bS_n
$$
The  differential is
given on the generators by
\Beq\label{LBk_infty}
\delta
\Ba{c}\resizebox{14mm}{!}{\begin{xy}
 <0mm,0mm>*{\circ};<0mm,0mm>*{}**@{},
 <-0.6mm,0.44mm>*{};<-8mm,5mm>*{}**@{-},
 <-0.4mm,0.7mm>*{};<-4.5mm,5mm>*{}**@{-},
 <0mm,0mm>*{};<-1mm,5mm>*{\ldots}**@{},
 <0.4mm,0.7mm>*{};<4.5mm,5mm>*{}**@{-},
 <0.6mm,0.44mm>*{};<8mm,5mm>*{}**@{-},
   <0mm,0mm>*{};<-8.5mm,5.5mm>*{^1}**@{},
   <0mm,0mm>*{};<-5mm,5.5mm>*{^2}**@{},
   <0mm,0mm>*{};<4.5mm,5.5mm>*{^{m\hspace{-0.5mm}-\hspace{-0.5mm}1}}**@{},
   <0mm,0mm>*{};<9.0mm,5.5mm>*{^m}**@{},
 <-0.6mm,-0.44mm>*{};<-8mm,-5mm>*{}**@{-},
 <-0.4mm,-0.7mm>*{};<-4.5mm,-5mm>*{}**@{-},
 <0mm,0mm>*{};<-1mm,-5mm>*{\ldots}**@{},
 <0.4mm,-0.7mm>*{};<4.5mm,-5mm>*{}**@{-},
 <0.6mm,-0.44mm>*{};<8mm,-5mm>*{}**@{-},
   <0mm,0mm>*{};<-8.5mm,-6.9mm>*{^1}**@{},
   <0mm,0mm>*{};<-5mm,-6.9mm>*{^2}**@{},
   <0mm,0mm>*{};<4.5mm,-6.9mm>*{^{n\hspace{-0.5mm}-\hspace{-0.5mm}1}}**@{},
   <0mm,0mm>*{};<9.0mm,-6.9mm>*{^n}**@{},
 \end{xy}}\Ea
\ \ = \ \
 \sum_{[1,\ldots,m]=I_1\sqcup I_2\atop
 {|I_1|\geq 0, |I_2|\geq 1}}
 \sum_{[1,\ldots,n]=J_1\sqcup J_2\atop
 {|J_1|\geq 1, |J_2|\geq 1}
}\hspace{0mm}
\pm
\Ba{c}\resizebox{22mm}{!}{ \begin{xy}
 <0mm,0mm>*{\circ};<0mm,0mm>*{}**@{},
 <-0.6mm,0.44mm>*{};<-8mm,5mm>*{}**@{-},
 <-0.4mm,0.7mm>*{};<-4.5mm,5mm>*{}**@{-},
 <0mm,0mm>*{};<0mm,5mm>*{\ldots}**@{},
 <0.4mm,0.7mm>*{};<4.5mm,5mm>*{}**@{-},
 <0.6mm,0.44mm>*{};<12.4mm,4.8mm>*{}**@{-},
     <0mm,0mm>*{};<-2mm,7mm>*{\overbrace{\ \ \ \ \ \ \ \ \ \ \ \ }}**@{},
     <0mm,0mm>*{};<-2mm,9mm>*{^{I_1}}**@{},
 <-0.6mm,-0.44mm>*{};<-8mm,-5mm>*{}**@{-},
 <-0.4mm,-0.7mm>*{};<-4.5mm,-5mm>*{}**@{-},
 <0mm,0mm>*{};<-1mm,-5mm>*{\ldots}**@{},
 <0.4mm,-0.7mm>*{};<4.5mm,-5mm>*{}**@{-},
 <0.6mm,-0.44mm>*{};<8mm,-5mm>*{}**@{-},
      <0mm,0mm>*{};<0mm,-7mm>*{\underbrace{\ \ \ \ \ \ \ \ \ \ \ \ \ \ \
      }}**@{},
      <0mm,0mm>*{};<0mm,-10.6mm>*{_{J_1}}**@{},
 <13mm,5mm>*{};<13mm,5mm>*{\circ}**@{},
 <12.6mm,5.44mm>*{};<5mm,10mm>*{}**@{-},
 <12.6mm,5.7mm>*{};<8.5mm,10mm>*{}**@{-},
 <13mm,5mm>*{};<13mm,10mm>*{\ldots}**@{},
 <13.4mm,5.7mm>*{};<16.5mm,10mm>*{}**@{-},
 <13.6mm,5.44mm>*{};<20mm,10mm>*{}**@{-},
      <13mm,5mm>*{};<13mm,12mm>*{\overbrace{\ \ \ \ \ \ \ \ \ \ \ \ \ \ }}**@{},
      <13mm,5mm>*{};<13mm,14mm>*{^{I_2}}**@{},
 <12.4mm,4.3mm>*{};<8mm,0mm>*{}**@{-},
 <12.6mm,4.3mm>*{};<12mm,0mm>*{\ldots}**@{},
 <13.4mm,4.5mm>*{};<16.5mm,0mm>*{}**@{-},
 <13.6mm,4.8mm>*{};<20mm,0mm>*{}**@{-},
     <13mm,5mm>*{};<14.3mm,-2mm>*{\underbrace{\ \ \ \ \ \ \ \ \ \ \ }}**@{},
     <13mm,5mm>*{};<14.3mm,-4.5mm>*{_{J_2}}**@{},
 \end{xy}}\Ea
\Eeq
where the signs on the r.h.s\ are uniquely fixed by the fact that they all equal to $+1$ for $d$ odd. On the other hand, the dg prop  $\HoLoBd$
 is generated
by the (skew)symmetric corollas of degree $1-d(m+n+2a-2)$,
$$
\Ba{c}\resizebox{16mm}{!}{\xy
(-9,-6)*{};
(0,0)*+{a}*\cir{}
**\dir{-};
(-5,-6)*{};
(0,0)*+{a}*\cir{}
**\dir{-};
(9,-6)*{};
(0,0)*+{a}*\cir{}
**\dir{-};
(5,-6)*{};
(0,0)*+{a}*\cir{}
**\dir{-};
(0,-6)*{\ldots};
(-10,-8)*{_1};
(-6,-8)*{_2};
(10,-8)*{_n};
(-9,6)*{};
(0,0)*+{a}*\cir{}
**\dir{-};
(-5,6)*{};
(0,0)*+{a}*\cir{}
**\dir{-};
(9,6)*{};
(0,0)*+{a}*\cir{}
**\dir{-};
(5,6)*{};
(0,0)*+{a}*\cir{}
**\dir{-};
(0,6)*{\ldots};
(-10,8)*{_1};
(-6,8)*{_2};
(10,8)*{_m};
\endxy}\Ea
=(-1)^{d(\sigma+\tau)}
\Ba{c}\resizebox{20mm}{!}{\xy
(-9,-6)*{};
(0,0)*+{a}*\cir{}
**\dir{-};
(-5,-6)*{};
(0,0)*+{a}*\cir{}
**\dir{-};
(9,-6)*{};
(0,0)*+{a}*\cir{}
**\dir{-};
(5,-6)*{};
(0,0)*+{a}*\cir{}
**\dir{-};
(0,-6)*{\ldots};
(-12,-8)*{_{\tau(1)}};
(-6,-8)*{_{\tau(2)}};
(12,-8)*{_{\tau(n)}};
(-9,6)*{};
(0,0)*+{a}*\cir{}
**\dir{-};
(-5,6)*{};
(0,0)*+{a}*\cir{}
**\dir{-};
(9,6)*{};
(0,0)*+{a}*\cir{}
**\dir{-};
(5,6)*{};
(0,0)*+{a}*\cir{}
**\dir{-};
(0,6)*{\ldots};
(-12,8)*{_{\sigma(1)}};
(-6,8)*{_{\sigma(2)}};
(12,8)*{_{\sigma(m)}};
\endxy}\Ea\ \ \ \forall \sigma\in \bS_m, \forall \tau\in \bS_n,
$$
where $m+n+ a\geq 3$, $m\geq 1$, $n\geq 1$, $a\geq 0$; the differential
 is given by
\Beq\label{2: d on Lie inv infty}
\delta
\Ba{c}\resizebox{16mm}{!}{\xy
(-9,-6)*{};
(0,0)*+{a}*\cir{}
**\dir{-};
(-5,-6)*{};
(0,0)*+{a}*\cir{}
**\dir{-};
(9,-6)*{};
(0,0)*+{a}*\cir{}
**\dir{-};
(5,-6)*{};
(0,0)*+{a}*\cir{}
**\dir{-};
(0,-6)*{\ldots};
(-10,-8)*{_1};
(-6,-8)*{_2};
(10,-8)*{_n};
(-9,6)*{};
(0,0)*+{a}*\cir{}
**\dir{-};
(-5,6)*{};
(0,0)*+{a}*\cir{}
**\dir{-};
(9,6)*{};
(0,0)*+{a}*\cir{}
**\dir{-};
(5,6)*{};
(0,0)*+{a}*\cir{}
**\dir{-};
(0,6)*{\ldots};
(-10,8)*{_1};
(-6,8)*{_2};
(10,8)*{_m};
\endxy}\Ea
=
\sum_{l\geq 1}\sum_{a=b+c+l-1}\sum_{[m]=I_1\sqcup I_2\atop
[n]=J_1\sqcup J_2} \pm
\Ba{c}
%
%
\Ba{c}\resizebox{21mm}{!}{\xy
(0,0)*+{b}*\cir{}="b",
(10,10)*+{c}*\cir{}="c",
%
(-9,6)*{}="1",
(-7,6)*{}="2",
(-2,6)*{}="3",
(-3.5,5)*{...},
(-4,-6)*{}="-1",
(-2,-6)*{}="-2",
(4,-6)*{}="-3",
(1,-5)*{...},
(0,-8)*{\underbrace{\ \ \ \ \ \ \ \ }},
(0,-11)*{_{J_1}},
(-6,8)*{\overbrace{ \ \ \ \ \ \ }},
(-6,11)*{_{I_1}},
(6,16)*{}="1'",
(8,16)*{}="2'",
(14,16)*{}="3'",
(11,15)*{...},
(11,6)*{}="-1'",
(16,6)*{}="-2'",
(18,6)*{}="-3'",
(13.5,6)*{...},
(15,4)*{\underbrace{\ \ \ \ \ \ \ }},
(15,1)*{_{J_2}},
(10,18)*{\overbrace{ \ \ \ \ \ \ \ \ }},
(10,21)*{_{I_2}},
%
(0,2)*-{};(8.0,10.0)*-{}
**\crv{(0,10)};
(0.5,1.8)*-{};(8.5,9.0)*-{}
**\crv{(0.4,7)};
(1.5,0.5)*-{};(9.1,8.5)*-{}
**\crv{(5,1)};
(1.7,0.0)*-{};(9.5,8.6)*-{}
**\crv{(6,-1)};
(5,5)*+{...};
\ar @{-} "b";"1" <0pt>
\ar @{-} "b";"2" <0pt>
\ar @{-} "b";"3" <0pt>
\ar @{-} "b";"-1" <0pt>
\ar @{-} "b";"-2" <0pt>
\ar @{-} "b";"-3" <0pt>
\ar @{-} "c";"1'" <0pt>
\ar @{-} "c";"2'" <0pt>
\ar @{-} "c";"3'" <0pt>
\ar @{-} "c";"-1'" <0pt>
\ar @{-} "c";"-2'" <0pt>
\ar @{-} "c";"-3'" <0pt>
\endxy}\Ea
\Ea
\Eeq
where the summation parameter $l$ counts the number of internal edges connecting the two vertices
on the r.h.s., and the signs are  fixed by the fact that they all equal to $+1$ for $d$ even. If $d=1$, it is called the prop of {\em strongly homotopy involutive Lie bialgebras}\, and is denoted by $\HoLoB$.

\sip

There is a canonical injection of props $\HoLBd\rar \HoLoBd$ sending a generator
$\HoLBd$ into the corresponding generator of $\HoLoBd$ with $a=0$.
\sip

Sometimes it is more suitable to work  with the degree shifted
 version $\HoLoBd\{d\}$ of the prop $\HoLoBd$ which is defined uniquely by the following property: a representation of $\HoLoBd\{d\}$ in a vector space $V$ is identical to the representation of $\HoLoBd$ in $V[d]$. The prop $\HoLoBd\{d\}$ is generated by the symmetric corollas of degree $1-2d(n+a-1)$,
$$
\Ba{c}\resizebox{16mm}{!}{\xy
(-9,-6)*{};
(0,0)*+{a}*\cir{}
**\dir{-};
(-5,-6)*{};
(0,0)*+{a}*\cir{}
**\dir{-};
(9,-6)*{};
(0,0)*+{a}*\cir{}
**\dir{-};
(5,-6)*{};
(0,0)*+{a}*\cir{}
**\dir{-};
(0,-6)*{\ldots};
(-10,-8)*{_1};
(-6,-8)*{_2};
(10,-8)*{_n};
(-9,6)*{};
(0,0)*+{a}*\cir{}
**\dir{-};
(-5,6)*{};
(0,0)*+{a}*\cir{}
**\dir{-};
(9,6)*{};
(0,0)*+{a}*\cir{}
**\dir{-};
(5,6)*{};
(0,0)*+{a}*\cir{}
**\dir{-};
(0,6)*{\ldots};
(-10,8)*{_1};
(-6,8)*{_2};
(10,8)*{_m};
\endxy}\Ea
=
\Ba{c}\resizebox{20mm}{!}{\xy
(-9,-6)*{};
(0,0)*+{a}*\cir{}
**\dir{-};
(-5,-6)*{};
(0,0)*+{a}*\cir{}
**\dir{-};
(9,-6)*{};
(0,0)*+{a}*\cir{}
**\dir{-};
(5,-6)*{};
(0,0)*+{a}*\cir{}
**\dir{-};
(0,-6)*{\ldots};
(-12,-8)*{_{\tau(1)}};
(-6,-8)*{_{\tau(2)}};
(12,-8)*{_{\tau(n)}};
(-9,6)*{};
(0,0)*+{a}*\cir{}
**\dir{-};
(-5,6)*{};
(0,0)*+{a}*\cir{}
**\dir{-};
(9,6)*{};
(0,0)*+{a}*\cir{}
**\dir{-};
(5,6)*{};
(0,0)*+{a}*\cir{}
**\dir{-};
(0,6)*{\ldots};
(-12,8)*{_{\sigma(1)}};
(-6,8)*{_{\sigma(2)}};
(12,8)*{_{\sigma(m)}};
\endxy}\Ea\ \ \ \forall \sigma\in \bS_m, \forall \tau\in \bS_n,
$$
The differential is given by the above formula with the slightly ambiguous symbol $\pm$ replaced by $+1$ (i.e.\ omitted); hence the sign rules become especially  simple in this case.

\subsection{\bf $\HoLoBd$-algebras as Maurer-Cartan elements \cite{CMW}}\label{2: MC as Holieb}
According to the general theory \cite{MV}, the set of $\HoLoBd\{d\}$-algebra structures in a dg
vector space $(V,\delta)$ can be identified with the set of Maurer-Cartan elements
  of a  graded Lie algebra,
\Beq\label{2: g_V diamond}
\fg_V^\diamond: = \Hom(V,V)[1]\ \oplus \ \Def\left(\HoLoBd\{d\} \stackrel{0}{\rar} \cE nd_V\right),
\Eeq
controlling deformations of the trivial morphism which sends all the generators
of $\HoLoBd\{d\}$ to zero in $\cE nd_V$. The summand $\Hom(V,V)[1]$
takes care about deformations of  the given differential $\delta$ in $V$.
Using the explicit description of the dg prop
$\HoLoBd\{d\}$ given at the end of the previous subsection, one can identify $\fg_V^\diamond$ as a $\Z$-graded vector space with
\Beqrn
\fg_V^\diamond
&=&\prod_{a\geq 0, m,n\geq 1} \Hom_{\bS_m\times \bS_n}\left(\Id_m\ot\Id_n[2d(n+a-1)-1],
\Hom\left(V^{\ot n}, V^{\ot  m}\right)\right)[-1]\\
&=& \prod_{a\geq0, m,n\geq 1}
\Hom\left((V[2d])^{\odot n}, V^{\odot  m}\right)[-2da+2d]
\Eeqrn
 Assume $V$ has a countable basis $(x_1, x_2, \ldots)$, and let $(p^1, p^2,...)$ stand for the associated set of dual generators of
$\Hom(V[2d],\K)$ (with $|p^i| + |x_i|=2d$), then
the degree shifted vector space $\fg_V^\diamond[-2d]$ can be identified
with the subspace of a graded commutative ring
$$
\fg_V^\diamond[-2d]\subset \K[[x^i,p_i,\hbar]]
$$
spanned by those formal power series $f(x,p,\hbar)$ which satisfy the  conditions
$$
f(x,p,\hbar)|_{x_i=0} =0, \ \ \ \ \  f(x,p,\hbar)|_{p^i=0} =0
$$
i.e. which belong to the maximal ideal generated by the products $x_ip^j$. Here $ \hbar$ is a formal  parameter\footnote{For a vector space $W$ we denote by $W[[\hbar]]$
 the vector space of formal power series in $\hbar$ with coefficients in $W$. For later use we denote by $\hbar^m W[[\hbar]]$ the subspace of $W[[\hbar]]$ spanned by series of the form
 $\hbar^m f$ for some $f\in W[[\hbar]]$.}
 of degree $2d$. The algebra $\K[[x^i,p_i,\hbar]]$ has a classical associative star product given explicitly (up to standard Koszul signs) as follows
$$
f *_\hbar g:=\sum_{k=0}^\infty \frac{\hbar^{k}}{k!}\sum_{i_1,\ldots, i_k}
\pm \frac{\p^k f}{\p p^{i_1}\cdots \p p^{i_k}}\frac{\p^ k g}{\p x_{i_1}\cdots
\p x_{i_k}}
$$
The Lie brackets in the degree shifted deformation complex $\fg_V^\diamond[-2d]$ are then given by \cite{DCTT}
$$
[f,g]=\frac{f *_\hbar g - (-1)^{|f||g|} g *_\hbar f}{\hbar}
$$
Hence $\HoLoBd\{d\}$-algebra structures in a graded vector space $V$ (with a countable basis) are in 1-1 correspondence with homogeneous formal power series $\Ga\in \fg_V^\diamond[-2d]$
of degree $1+2d$ satisfying the equation
\Beq\label{2: equation for Gamma-h}
\Ga *_\hbar \Ga=\sum_{k=1}^\infty \frac{\hbar^{k-1}}{k!}\sum_{i_1,\ldots, i_k}
\pm \frac{\p^ k \Ga}{\p p^{i_1}\cdots p^{i_k}}\frac{\p^ k \Ga}{\p x_{i_1}\cdots
\p x_{i_k}}=0,
\Eeq
This compact description of all higher homotopy involutive Lie bialgebra operations in $V$ is quite useful in making an explicit link between
$\HoLoBd$ algebras and {\em commutative $BV$-algebras}\, which is outlined next.

\subsection{Commutative Batalin-Vilkovisky $d$-algebras}\label{2: Comm BV}
A {\em commutative Batalin-Vilkovisky $d$-algebra}\, or, shortly, a $\HocBVd$-{\em algebra}\, is, by definition \cite{Kr}, a differential graded commutative algebra  $(V,\delta)$
equipped with a countable collections of homogeneous linear maps, $\{\Delta_a: V\rar V,\
    |\Delta_a|=1-2da\}_{a\geq 1}$, such that each operator $\Delta_a$ is of order $\leq a+1$ (with respect to the given multiplication) and the equations,
    \Beq\label{5: BV_comm equation for Delta}
    \sum_{a=0}^n \Delta_a \circ \Delta_{n-a}=0, \ \ \ \text{with}\ \Delta_0:=-\delta
    \Eeq
hold for any $n\in \N$. These equations are equivalent to one equation,
    $$
    \Delta_\hbar^2=0
    $$
for the formal power series of operators
$$
\Delta_\hbar:=\sum_{a=0}^\infty \hbar^a \Delta_a,
$$
where the formal power variable $\hbar$ is assigned degree $2d$.
Let us denote by
 $\HocBVd$ the dg operad governing commutative $BV$ $d$-algebras. This operad
is the quotient
of the free operad generated by one binary operation of degree zero
$
\begin{xy}
 <0mm,0.66mm>*{};<0mm,3mm>*{}**@{-},
 <0.39mm,-0.39mm>*{};<2.2mm,-2.2mm>*{}**@{-},
 <-0.35mm,-0.35mm>*{};<-2.2mm,-2.2mm>*{}**@{-},
 <0mm,0mm>*{\circ};<0mm,0mm>*{}**@{},
   <0.39mm,-0.39mm>*{};<2.9mm,-4mm>*{^2}**@{},
   <-0.35mm,-0.35mm>*{};<-2.8mm,-4mm>*{^1}**@{},
\end{xy}=
\begin{xy}
 <0mm,0.66mm>*{};<0mm,3mm>*{}**@{-},
 <0.39mm,-0.39mm>*{};<2.2mm,-2.2mm>*{}**@{-},
 <-0.35mm,-0.35mm>*{};<-2.2mm,-2.2mm>*{}**@{-},
 <0mm,0mm>*{\circ};<0mm,0mm>*{}**@{},
   <0.39mm,-0.39mm>*{};<2.9mm,-4mm>*{^1}**@{},
   <-0.35mm,-0.35mm>*{};<-2.8mm,-4mm>*{^2}**@{},
\end{xy}
$ (standing for graded commutative multiplication)
and a countable family of unary operations
$
\left\{ \Ba{c}\resizebox{3.1mm}{!}{  \xy
(0,5)*{};
(0,0)*+{_a}*\cir{}
**\dir{-};
(0,-5)*{};
(0,0)*+{_a}*\cir{}
**\dir{-};
\endxy}\Ea \right\}_{a\geq 1}
$
of homological degree $1-2da$ modulo the ideal $I$ generated by
the standard associativity relation
for  $\begin{xy}
 <0mm,0.66mm>*{};<0mm,3mm>*{}**@{-},
 <0.39mm,-0.39mm>*{};<2.2mm,-2.2mm>*{}**@{-},
 <-0.35mm,-0.35mm>*{};<-2.2mm,-2.2mm>*{}**@{-},
 <0mm,0mm>*{\circ};<0mm,0mm>*{}**@{},
\end{xy}$ and the compatibility relations involving the latter and the unary operations
which assure that each unary operation $\Ba{c}\resizebox{2.9mm}{!}{\xy
(0,5)*{};
(0,0)*+{_a}*\cir{}
**\dir{-};
(0,-5)*{};
(0,0)*+{_a}*\cir{}
**\dir{-};
\endxy}\Ea$  is of order $\leq a+1$ with
respect to the multiplication.
The differential $\delta$ in $\BV_\infty^{com}$ is given by
\Beq\label{5: diff in BV comm}
\delta\begin{xy}
 <0mm,0.66mm>*{};<0mm,3mm>*{}**@{-},
 <0.39mm,-0.39mm>*{};<2.2mm,-2.2mm>*{}**@{-},
 <-0.35mm,-0.35mm>*{};<-2.2mm,-2.2mm>*{}**@{-},
 <0mm,0mm>*{\circ};<0mm,0mm>*{}**@{},
   <0.39mm,-0.39mm>*{};<2.9mm,-4mm>*{^2}**@{},
   <-0.35mm,-0.35mm>*{};<-2.8mm,-4mm>*{^1}**@{},
\end{xy}=0, \ \ \ \
\delta\ \Ba{c}\resizebox{4mm}{!}{ \xy
(0,5)*{};
(0,0)*+{_a}*\cir{}
**\dir{-};
(0,-5)*{};
(0,0)*+{_a}*\cir{}
**\dir{-};
\endxy}\Ea: =\sum_{a=b+c\atop b,c\geq 1}\Ba{c} \resizebox{4mm}{!}{
\xy
(0,6,3)*{};
(0,0)*+{_c}*\cir{}
**\dir{-};
(0,-5)*{};
(0,0)*+{_c}*\cir{}
**\dir{-};
(0,13)*{};
(0,8)*+{_b}*\cir{}
**\dir{-};
\endxy}
\Ea
\Eeq
Let $J$ be the differential closure of an ideal in $\HocBVd$ generated by operations $\Ba{c}\resizebox{2.9mm}{!}{\xy
(0,5)*{};
(0,0)*+{_a}*\cir{}
**\dir{-};
(0,-5)*{};
(0,0)*+{_a}*\cir{}
**\dir{-};
\endxy}\Ea$ with $a\geq 2$. The quotient $\HocBVd/J$ is precisely of the operad of (degree shifted) Batalin-Vilkovisky algebras $\BV_d$.
It was proven in \cite{CMW} that the canonical projection
$
\HocBVd \lon \BV_d$
is quasi-isomorphism of
operads.

\subsection{From $\HoLoBd$-algebras to  $\HocBVd$-algebras and back}\label{2: Com BV to Holieb}
 Recall that a $\HoLoBd\{d\}$-algebra structure in a graded vector space
$V$ (i.e.\ a $\HoLoBd$ structure in $V[d]$) can be identified with a degree $1+2d$ element $\Ga$ in $\fg_V^\diamond[-2d]\subset \K[[p^i, x_i,\hbar]]$
satisfying equation  (\ref{2: equation for Gamma-h}). Out of this datum
one creates a
$\HocBVd$-algebra structure on $\widehat{\odot^\bu}V$ (the completed symmetric tensor algebra on $V$), i.e.\ a representation
$$
\rho: \HocBVd \lon \cE nd_{\widehat{\odot^\bu}(V)},
$$
which is given explicitly as follows \cite{CMW},
$$
\left\{
\Ba{l}
\rho\left(\begin{xy}
 <0mm,0.66mm>*{};<0mm,3mm>*{}**@{-},
 <0.39mm,-0.39mm>*{};<2.2mm,-2.2mm>*{}**@{-},
 <-0.35mm,-0.35mm>*{};<-2.2mm,-2.2mm>*{}**@{-},
 <0mm,0mm>*{\circ};<0mm,0mm>*{}**@{},
\end{xy}\right):=\mathrm{the\ standard\ multiplication\ in}\
\widehat{{\odot}^\bu}(V)\\
\Delta_a:=\rho\left(\Ba{c}\resizebox{4mm}{!}{ \xy
(0,5)*{};
(0,0)*+{_a}*\cir{}
**\dir{-};
(0,-5)*{};
(0,0)*+{_a}*\cir{}
**\dir{-};
\endxy}\Ea\right):= \displaystyle\sum_{a+1=k+l} \frac{1}{k!l!}\sum_{i_1,\ldots,i_{l}} \frac{\p^{a+1} \Ga}{\p^k\hbar \p p^{i_1}\cdots
\p p^{i_{l}}}|_{\hbar=p^i=0}\frac{\p^{l} }{\p x_{i_1}\cdots
\p x_{i_{l}}}
\Ea
\right.
$$
This explicit correspondence can be equivalently understood as  a morphism of dg operads
$$
F: \HocBVd \lon \f (\HoLoBd\{d\})
$$
given explicitly on generators as follows,
$$
F\left(
\begin{xy}
 <0mm,0.66mm>*{};<0mm,3mm>*{}**@{-},
 <0.39mm,-0.39mm>*{};<2.2mm,-2.2mm>*{}**@{-},
 <-0.35mm,-0.35mm>*{};<-2.2mm,-2.2mm>*{}**@{-},
 <0mm,0mm>*{\circ};<0mm,0mm>*{}**@{},
   <0.39mm,-0.39mm>*{};<2.9mm,-4mm>*{^2}**@{},
   <-0.35mm,-0.35mm>*{};<-2.8mm,-4mm>*{^1}**@{},
\end{xy}\right):=
\Ba{c}\resizebox{9mm}{!}{
\xy
(0,4)*+{_{\ \, }}*\frm{o};
(-3,-3)*+{_{_1}}*\frm{o};
(3,-3)*+{_{_2}}*\frm{o};
\endxy}\Ea
$$
\Beq\label{2: map F}
F\left(\Ba{c}\resizebox{4mm}{!}{\xy
(0,5)*{};
(0,0)*+{_a}*\cir{}
**\dir{-};
(0,-5)*{};
(0,0)*+{_a}*\cir{}
**\dir{-};
\endxy}\Ea\right)
=\sum_{m\geq 1}\sum_{a+1=p+q}
\Ba{c}\resizebox{8mm}{!}{ \xy
(0,6)*{_m},
(0,-6)*{^l},
(0,4)*{_{...}},
(0,-4)*{^{...}},
 (0,-10)*+{\ }*\frm{o}="B";
  (0,0)*+{_{_k}}*\frm{o}="C";
  (0,10)*+{\ }*\frm{o}="X";
"C"; "X" **\crv{(-8,3)}; 
 "C"; "X" **\crv{(-5,3)};
  "C"; "X" **\crv{(5,3)};
   "C"; "X" **\crv{(8,3)};
   "C"; "B" **\crv{(-8,-3)}; 
 "C"; "B" **\crv{(-5,-3)};
  "C"; "B" **\crv{(5,-3)};
   "C"; "B" **\crv{(8,-3)};
 \endxy}
 \Ea
\Eeq
where
$$
\Ba{rccc}
\f: & \text{Category of props} & \lon & \text{Category of operads}\\
    &   \cP & \lon & \f(\cP)
    \Ea
$$
is the {\em polydifferential functor}\, introduced in \cite{MW} (we refer to \S 5.1 of \cite{MW} or \S 2.2 of \cite{MW2}
for full details explaining, in particular, the symbols on the r.h.s.\ of the above formula; these sections can be read independently of the rest of both papers). Its main defining property is that, given any representation $\rho:\cP \rar \cE nd_V$ of a prop $\cP$ in a graded vector space $X$, there is an associated representation $\f(\rho): \f(\cP)\rar \cE nd_{\widehat{\odot^\bu} X}$ of the operad $\f(\cP)$ in the completed free graded commutative algebra $\widehat{\odot^\bu}X$ such that
elements of $\cP$ acts on $\widehat{\odot^\bu}(X)$ as polydifferential operators. The symbol on the r.h.s.\ of (\ref{2: map F}) is precisely the polydifferential operator corresponding to the generator $\Ba{c}\resizebox{13mm}{!}{\xy
(-9,-6)*{};
(0,0)*+{_k}*\cir{}
**\dir{-};
(-5,-6)*{};
(0,0)*+{_k}*\cir{}
**\dir{-};
(9,-6)*{};
(0,0)*+{_k}*\cir{}
**\dir{-};
(5,-6)*{};
(0,0)*+{_k}*\cir{}
**\dir{-};
(0,-6)*{\ldots};
(-10,-8)*{_1};
(-6,-8)*{_2};
(10,-8)*{_l};
(-9,6)*{};
(0,0)*+{_k}*\cir{}
**\dir{-};
(-5,6)*{};
(0,0)*+{_k}*\cir{}
**\dir{-};
(9,6)*{};
(0,0)*+{_k}*\cir{}
**\dir{-};
(5,6)*{};
(0,0)*+{_k}*\cir{}
**\dir{-};
(0,6)*{\ldots};
(-10,8)*{_1};
(-6,8)*{_2};
(10,8)*{_m};
\endxy}\Ea$ of $\HoLoBd\{d\}$. Reversely, given a representation
$$
\rho: \HocBVd \lon \cE nd_{\widehat{\odot^\bu}(V)},
$$
such that
$$
\rho\left(\begin{xy}
 <0mm,0.66mm>*{};<0mm,3mm>*{}**@{-},
 <0.39mm,-0.39mm>*{};<2.2mm,-2.2mm>*{}**@{-},
 <-0.35mm,-0.35mm>*{};<-2.2mm,-2.2mm>*{}**@{-},
 <0mm,0mm>*{\circ};<0mm,0mm>*{}**@{},
\end{xy}\right):=\mathrm{the\ standard\ multiplication\ in}\
\widehat{{\odot}^\bu}(V)
$$
it was proven in \cite{DCTT} that it factors through the composition
$$
\rho:  \HocBVd \stackrel{F}{\lon} \f (\HoLoBd\{d\}) \stackrel{\f(\rho')}{\lon}
\f(\cE nd_V)= \cE nd_{\widehat{\odot^\bu} V}
$$
for some representation $\rho': \HoLoB\{d\}\rar \cE nd_V$. Put another way, one can read all the higher homotopy involutive Lie bialgebra operations $\rho'\left(\Ba{c}\resizebox{13mm}{!}{\xy
(-9,-6)*{};
(0,0)*+{_k}*\cir{}
**\dir{-};
(-5,-6)*{};
(0,0)*+{_k}*\cir{}
**\dir{-};
(9,-6)*{};
(0,0)*+{_k}*\cir{}
**\dir{-};
(5,-6)*{};
(0,0)*+{_k}*\cir{}
**\dir{-};
(0,-6)*{\ldots};
(-10,-8)*{_1};
(-6,-8)*{_2};
(10,-8)*{_l};
(-9,6)*{};
(0,0)*+{_k}*\cir{}
**\dir{-};
(-5,6)*{};
(0,0)*+{_k}*\cir{}
**\dir{-};
(9,6)*{};
(0,0)*+{_k}*\cir{}
**\dir{-};
(5,6)*{};
(0,0)*+{_k}*\cir{}
**\dir{-};
(0,6)*{\ldots};
(-10,8)*{_1};
(-6,8)*{_2};
(10,8)*{_m};
\endxy}\Ea\right)$ in $V$ from the explicit representation of
$$
\rho\left(
\Ba{c}\resizebox{4mm}{!}{\xy
(0,5)*{};
(0,0)*+{_a}*\cir{}
**\dir{-};
(0,-5)*{};
(0,0)*+{_a}*\cir{}
**\dir{-};
\endxy}\Ea
\right)=
\sum_{m\geq 1}\frac{1}{m!}\sum_{a+1=k+l}\sum_{i_\bu,j_\bu} C^{(k)\, j_1\ldots j_m}_{\ \ \ i_1\ldots i_l} x_{j_1}\ldots x_{j_m} \frac{\p^{l} }{\p x_{i_1}\cdots
\p x_{i_{l}}}
$$
of the generator $\Ba{c}\resizebox{3mm}{!}{\xy
(0,5)*{};
(0,0)*+{_a}*\cir{}
**\dir{-};
(0,-5)*{};
(0,0)*+{_a}*\cir{}
**\dir{-};
\endxy}\Ea$
as a differential operator of order $\leq a+1$ on $\widehat{\odot^\bu} V\simeq
\K[[x_i]]$: the linear map
$$
\mu_{\ \ l}^{(k) m}:=\rho'\left(\Ba{c}\resizebox{13mm}{!}{\xy
(-9,-6)*{};
(0,0)*+{_k}*\cir{}
**\dir{-};
(-5,-6)*{};
(0,0)*+{_k}*\cir{}
**\dir{-};
(9,-6)*{};
(0,0)*+{_k}*\cir{}
**\dir{-};
(5,-6)*{};
(0,0)*+{_k}*\cir{}
**\dir{-};
(0,-6)*{\ldots};
(-10,-8)*{_1};
(-6,-8)*{_2};
(10,-8)*{_l};
(-9,6)*{};
(0,0)*+{_k}*\cir{}
**\dir{-};
(-5,6)*{};
(0,0)*+{_k}*\cir{}
**\dir{-};
(9,6)*{};
(0,0)*+{_k}*\cir{}
**\dir{-};
(5,6)*{};
(0,0)*+{_k}*\cir{}
**\dir{-};
(0,6)*{\ldots};
(-10,8)*{_1};
(-6,8)*{_2};
(10,8)*{_m};
\endxy}\Ea\right): \ot^m V \rar \ot^n V
$$
is given in the basis $\{x_i\}$ by (modulo the standard Koszul signs) by
$$
\mu_{\ \ \ l}^{(k) m}(x_{i_1}\ot \ldots \ot x_{i_l})=\sum_{j_\bu} C^{(k)\, j_1\ldots j_m}_{\ \ \ i_1\ldots i_l} x_{j_1}\ot \ldots\ot  x_{j_m}
$$
We shall use this one-to-one correspondence heavily in the proof of the Main Theorem {\ref{4: Main Th}} below.

\subsection{A basic example of an involutive Lie bialgebra}\label{2: subsection on cyclic words} Let $W$ be a graded vector space equipped with a
linear map
$\Theta:  \odot^2(W[d]) \rar \K[1+d]$ for some $d\in \Z$. This map is the
same as a degree $1-d$ pairing
$\Theta:  W\ot W \rar \K[1-d]$
satisfying the (skew) symmetry condition,
\Beq\label{2: skewsym of Theta}
\Theta(w_1,w_2)=(-1)^{d+|w_1||w_2|}\Theta(w_2,w_1), \ \ \forall\ w_1,w_2\in W.
\Eeq
A symplectic structure on $W$ corresponds to the case $d=1$ and $\Theta$ non-degenerate.

\sip

 The associated
vector space of ``cyclic words in $W$",
$$
Cyc^\bu(W):=\sum_{n\geq 0} (W^{\ot n})^{\Z_n},
$$
admits a $\LoBd$-structure given by the following well-known
formulae for the Lie bracket and cobracket,
$$
[(w_1\ot...\ot w_n)_{\Z_n}, (w_1'\ot ...\ot w'_m)^{\Z_n}]:=\hspace{100mm}
$$
$$
\hspace{20mm} \sum_{i=1}^n\sum_{j=1}^m
 \pm
\Theta(w_i,w'_j) (w_1\ot ...\ot  w_{i-1}\ot w'_{j+1}\ot ... \ot w'_m\ot w'_1\ot ... \ot w'_{j-1}\ot w_{i+1}\ot\ldots\ot w_n)^{\Z_{n+m-2}}
$$
$$
\vartriangle (w_1\ot\ldots\ot w_n)_{\Z_n}:=\sum_{i\neq j}
\pm\Theta(w_i,w_j)(w_{i+1}\ot ...\ot w_{j-1})^{\Z_{j-i-1}}\bigotimes
(w_{j+1}\ot ...\ot w_{i-1})^{\Z_{n-j+i-1}}    \hspace{15mm}
$$
where $\pm$ stands for the standard Koszul sign. A very short (and pictorial) proof of this claim can be found in \cite{MW}.
 Note that the vector space
$Cyc^\bu(W)$ is naturally {\em weight}\,-graded
$$
Cyc^\bu(W)= \bigoplus_{n\geq 0} Cyc^n(W),\ \ \ \ Cyc^n(W):=(\ot^n W)^{\Z_n}
$$
by the length of cyclic words,
and both operations $\Delta$ and $[\ ,\ ]$ have weight-degree $-2$ with respect to this weight-grading (which should not be confused with the {\em homological}\, grading).

\subsubsection{\bf A special case: Schedler's necklace Lie bialgebra}\label{2: Sch liab}
A special case of the above construction for $d=1$ gives us Schedler's necklace Lie bialgebra structure \cite{Sch} associated with the quiver
(\ref{1: quiver for W_N}). Consider a set of $N$ formal letters
$$
\{x_1,\ldots, x_n\}
$$
and denote by $W_N$ their linear span over a field $\K$.
We shall make
$Cyc(W)$ into a weight-degree $-1$ (not $-2$ as in the above example!) involutive Lie
bialgebra using the following ``doubling" trick.

\sip

Consider two copies $W_N^{(1)}, W_N^{(2)}$ of $W_N$ and equip their direct sum
$$
\hat{W}_N:= W_N^{(1)}\oplus W_N^{(2)}
$$
with the unique symplectic structure $\theta: \wedge^2 \hat{W}\rar \K$
 making the basis $\{x^{(1)}_\al, x^{(2)}_\be\}_{1\leq \al,\be\leq N}$  a Darboux one,
 $$
 \Theta(x^{(1)}_\al,x^{(2)}_\be)=-\Theta(x^{(2)}_\be,x^{(1)}_\al)=\delta_\be^\al,
 \ \ \Theta(x^{(1)}_\al,x^{(1)}_\be)=0, \ \ \ \ \Theta(x^{(2)}_\al,x^{(2)}_\be)=0.
 $$
Then the above formulae for $[\ ,\ ]$ and $\Delta$ make the space $Cyc(\hat{W}_N)$
into an involutive Lie bialgebra. It is easy to see that the subspace
$$
Cyc(\hat{W}_N)^{(12)}\subset Cyc(\hat{W}_N)
$$
spanned by cyclic words of the form
$$
\left(x_{\al_1}^{(1)}\ot x_{\al_1}^{(2)}\ot x_{\al_2}^{(1)}\ot x_{\al_2}^{(2)}\ot \ldots \ot
x_{\al_n}^{(1)}\ot x_{\al_n}^{(2)}\right)^{\Z_{2n}}
$$
is closed with respect to the above Lie bracket and co-bracket. The canonical isomorphism
$$
\Ba{ccc}
Cyc(W_N) & \lon & Cyc(\hat{W}_N)^{(12)}\\
\left(x_{\al_1}\ot x_{\al_2}\ot  \ldots \ot
x_{\al_n}\right)^{\Z_{n}} & \lon & \left(x_{\al_1}^{(1)}\ot x_{\al_1}^{(2)}\ot x_{\al_2}^{(1)}\ot x_{\al_2}^{(2)}\ot \ldots \ot
x_{\al_n}^{(1)}\ot x_{\al_n}^{(2)}\right)^{\Z_{2n}}
\Ea
$$
makes $Cyc(W_N)$ into an involutive Lie bialgebra with the
Lie bracket identical to the one introduced earlier by Schedler in \cite{Sch}, $[\ ,\ ]=[\ ,\ ]^{S}$, but with the Lie cobracket slightly different,
$$
\Delta \left(e_{\al_1}\ot  \ldots \ot
x_{\al_n}\right)^{\Z_{n}}= \Delta^{S} \left(x_{\al_1}\ot  \ldots \ot
x_{\al_n}\right)^{\Z_{n}} + \sum_{i=1}^n
1\wedge \left(x_{\al_{1}}\ot ...\ot x_{\al_{i-1}}\ot x_{\al_{i+1}}\ot \ldots \ot x_{\al_{n}}\right)^{\Z_{n-1}}.
$$
This  purely combinatorial Lie bialgebra structure
 on $Cyc(W_N)$ has a  beautiful geometric  interpretation ---
it is {\em isomorphic}\, \cite{AKKN, AN, Ma}
to the Goldman-Turaev Lie bialgebra
structure on the space of free loops in $\Sigma_{0,N+1}$, the two dimensional sphere  with
$n+1$ non-intersecting open disks removed.

\bip

\bip

\bip


{\large
\section{\bf A prop of ribbon hypergraphs}
}

\subsection{Ribbon hypergraphs}
A {\em ribbon hypergraph}\, $\Ga$ is, by definition, a triple
$(E(\Ga), \sigma_1, \sigma_0)$ consisting of a finite set  $E(\Ga)$
  {\em of edges}\, and two arbitrary bijections (``permutations") $\sigma_0,\sigma_1: E(\Ga) \rar E(\Ga)$.  The
orbits
$$
V(\Ga):=E(\Ga)/\sigma_0
$$
or, equivalently, the cycles of the permutation $\sigma_0$ are called the {\em vertices}\, of $\Ga$  while the orbits
$$
H(\Ga):=E(\Ga)/\sigma_1
$$
 of the permutation $\sigma_1$ are called 
 {\em hyperedges} of $\Ga$ (cf.\ \cite{LZ}).
 Let $p_\circ: E(\Ga)\rar V(\Ga)$ and $p_*: E(\Ga)\rar H(\Ga)$ be  canonical projections. For any  vertex $v\in V(\Ga)$ and any hyperedge $h\in H(\Ga)$ the associated sets of edges
 $$
 p_\circ^{-1}(v)=\{e_{i_1},\ldots, e_{i_k}\}, \ \ \ \  p_*^{-1}(h)=\{e_{j_1},\ldots, e_{j_l}\}, \ \ \ k,l\in \N
 $$
come equipped with induced cyclic orderings and hence can be represented pictorially as planar corollas,
$$
v \  \Leftrightarrow\ \
 \xy
(6.1,0)*{_{_{e_{i_1}}}},
(4.6,-4.5)*{_{_{e_{i_2}}}},
(0,-5)*{_{_{...}}},
 (0,0)*{\circ}="a",
(4,0)*{}="b1",
(3,-3)*{}="b2",
(0,-4)*{}="b3",
(-3,-3)*{}="b4",
(-4,0)*{}="b5",
(-3,3)*{}="b6",
(0,4)*{}="b7",
(3,3)*{}="b8",
\ar @{-} "a";"b1" <0pt>
\ar @{-} "a";"b2" <0pt>
\ar @{-} "a";"b3" <0pt>
\ar @{-} "a";"b4" <0pt>
\ar @{-} "a";"b5" <0pt>
\ar @{-} "a";"b6" <0pt>
\ar @{-} "a";"b7" <0pt>
\ar @{-} "a";"b8" <0pt>
\endxy\ \ , \ \ \
h \  \Leftrightarrow\ \
 \xy
(6.1,0)*{_{_{e_{j_1}}}},
(4.6,-4.5)*{_{_{e_{j_2}}}},
(0,-5)*{_{_{...}}},
 (0,0)*{\ast}="a",
(4,0)*{}="b1",
(3,-3)*{}="b2",
(0,-4)*{}="b3",
(-3,-3)*{}="b4",
(-4,0)*{}="b5",
(-3,3)*{}="b6",
(0,4)*{}="b7",
(3,3)*{}="b8",
\ar @{-} "a";"b1" <0pt>
\ar @{-} "a";"b2" <0pt>
\ar @{-} "a";"b3" <0pt>
\ar @{-} "a";"b4" <0pt>
\ar @{-} "a";"b5" <0pt>
\ar @{-} "a";"b6" <0pt>
\ar @{-} "a";"b7" <0pt>
\ar @{-} "a";"b8" <0pt>
\endxy
$$
Each edge $e\in E(\Ga)$ belongs to precisely one vertex, $e\in p_0^{-1}(v)$ for some $v\in V(\Ga)$,
and precisely one hyperedge, $e\in p_*^{-1}(h)$ for some $h\in H(\Ga)$. Hence we can glue the vertex $v$ to the corresponding hyperedge $h$ along the the common edge(s) to get a pictorial ($CW$ complex like) representation of a ribbon hypergraph $\Ga$ as an ordinary ribbon graph whose vertices are bicoloured, e.g.
$$
\Ga_1=\xy
 (0,0)*{\ast}="a",
(5,0)*{\circ}="b1",
(-3,4)*{\circ}="b2",
(-3,-4)*{\circ}="b3",
\ar @{-} "a";"b1" <0pt>
\ar @{-} "a";"b2" <0pt>
\ar @{-} "a";"b3" <0pt>
\endxy \ \ \ \ , \ \ \ \
\Ga_2=
\xy
(-4,0)*{\circ}="1";
(4,0)*{\star}="2";
"1";"2" **\crv{(-4,4) & (4,4)};
"1";"2" **\crv{(-4,-4) & (4,-4)};
\ar @{-} "1";"2" <0pt>
\endxy
\ \ \ \ , \ \ \ \
\Ga_3=
\xy
(-4,0)*{\circ}="1";
(4,0)*{\ast}="2";
(-10,0)*{\ast}="3";
"1";"2" **\crv{(-4,4) & (4,4)};
"1";"2" **\crv{(-4,-4) & (4,-4)};
\ar @{-} "1";"3" <0pt>
\endxy
$$
The vertices of $\Ga$ get represented pictorially as white vertices while hyperedges as asterisk vertices; sometimes we call a hyperedge an {\em asterisk vertex}\, when commenting some pictures. Note that each edge $e\in E(\Ga)$ connects precisely one white vertex to precisely one asterisk vertex; such bicoloured ribbon graphs are called {\em hypermaps}\, in \cite{LZ}.

\sip

Note that the set of edges $E(\Ga)$ admits two decompositions into disjoint subsets,
$$
E(\Ga)=\coprod_{v\in V(\Ga)} p_\circ^{-1}(v), \ \ \ \ E(\Ga)=\coprod_{h\in H(\Ga)} p_\ast^{-1}(h)
$$
with each subset having an induced cyclic ordering.

\sip

The orbits
 of the permutation $\sigma_\infty:= \sigma_0^{-1}\circ \sigma_1$
are called {\em boundaries}\, of the ribbon hypergraph $\Ga$; the set of boundaries is denoted
by $B(\Ga)$.
 For example, in the case of the above ribbon hypergraphs we
have $\# V(\Ga_1)=3$, $\# H(\Ga_1)=1$, $\# B(\Ga_1)=1$,  $\# E(\Ga_1)=3$,
 $\# V(\Ga_2)=1$, $\# H(\Ga_2)=1$, $\# B(\Ga_2)=3$,
 $\# E(\Ga_2)=3$  and $\# V(\Ga_3)=1$, $\# H(\Ga_3)=2$, $\# B(\Ga_3)=2$,  $\# E(\Ga_3)=3$. For a vertex $v\in V(\Ga)$ (resp, a hyperedge $h\in H(\Ga)$) we denote its {\em valency}\, by $|v|:=\# p_\circ^{-1}(v)$ (resp.,
 $|h|:=\# p_*^{-1}(h)$).

\sip

A ribbon hypergraph with each hyperedge having valency $2$ is called a {\em ribbon graph}. Ribbon graphs are depicted pictorially with asterisk vertices omitted as they contain no extra information,
$$
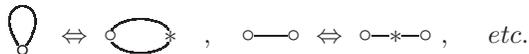
\xy
(0,-2)*{\circ}="A";
(0,-2)*{\circ}="B";
"A"; "B" **\crv{(6,6) & (-6,6)};
\endxy\ \Leftrightarrow \ \xy
(-4,0)*{\circ}="1";
(4,0)*{\ast}="2";
"1";"2" **\crv{(-4,3) & (4,3)};
"1";"2" **\crv{(-4,-3) & (4,-3)};
\endxy
 \ \ \  , \ \ \
 \xy
 (0,0)*{\circ}="a",
(6,0)*{\circ}="b",
\ar @{-} "a";"b" <0pt>
\endxy \ \Leftrightarrow \
\xy
 (0,0)*{\ast}="a",
(4,0)*{\circ}="b",
(-4,0)*{\circ}="c",
\ar @{-} "a";"b" <0pt>
\ar @{-} "a";"c" <0pt>
\endxy
\ , \ \ \ \ etc.
$$

\subsection{On geometric interpretation of ribbon hypergraphs}
Every connected ribbon graph $\Ga$
 can be interpreted geometrically as a topological 2-dimensional surface  with
 $\# B(\Ga)$ boundary circles and $\# V(\Ga)$ punctures which is obtained from its $CW$-complex
 realization by  thickening its every vertex into a closed disk punctured in the center and then  thickening its
every edge $e\in E(\Ga)$ into a 2-dimensional strip. For example
 (cf.\ \cite{MW}),
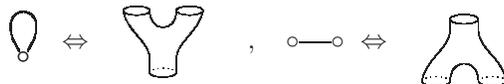
 $$
 \xy
(0,-2)*{\circ}="A";
(0,-2)*{\circ}="B";
"A"; "B" **\crv{(6,6) & (-6,6)};
\endxy\ \Leftrightarrow\
\Ba{c}\resizebox{12mm}{!}{\mbox{
\xy
(0,-1)*\ellipse(3,1){.};
(0,-1)*\ellipse(3,1)__,=:a(-180){-};
(-3,6)*\ellipse(3,1){-};
(3,6)*\ellipse(3,1){-};
(-3,12)*{}="1"; 
(3,12)*{}="2"; 
(-9,12)*{}="A2";
(9,12)*{}="B2"; 
"1";"2" **\crv{(-3,7) & (3,7)};
(-3,-2)*{}="A";
(3,-2)*{}="B"; 
(-3,1)*{}="A1";
(3,1)*{}="B1"; 
"A";"A1" **\dir{-};
"B";"B1" **\dir{-}; 
"B2";"B1" **\crv{(8,7) & (3,5)};
"A2";"A1" **\crv{(-8,7) & (-3,5)};
\endxy}
}\Ea
\ \ \ , \ \ \
 \xy
 (0,0)*{\circ}="a",
(6,0)*{\circ}="b",
\ar @{-} "a";"b" <0pt>
\endxy \ \Leftrightarrow \
\Ba{c}\resizebox{12mm}{!}{\mbox{
 \xy
(0,0.6)*\ellipse(3,1){-};
(-3,-6)*\ellipse(3,1){.};
(-3,-6)*\ellipse(3,1)__,=:a(-180){-};
(3,-6)*\ellipse(3,1){.};
(3,-6)*\ellipse(3,1)__,=:a(-180){-};
(-3,-12)*{}="1";
(3,-12)*{}="2";
(-9,-12)*{}="A2";
(9,-12)*{}="B2";
"1";"2" **\crv{(-3,-7) & (3,-7)};
(-3,1)*{}="A";
(3,1)*{}="B";
(-3,-1)*{}="A1";
(3,-1)*{}="B1";
"A";"A1" **\dir{-};
"B";"B1" **\dir{-};
"B2";"B1" **\crv{(8,-7) & (3,-5)};
"A2";"A1" **\crv{(-8,-7) & (-3,-5)};
\endxy}}\Ea
$$
with punctures represented as  bottom ``in-circles", and boundaries as top
``out-circles".

\sip

A ribbon hypergraph $\Ga$
is often used \cite{LZ} to encode combinatorially a Belyi map, that is, a ramified covering $f: X\rar \CP^1$ of the sphere whose ramification locus is contained in the set $\{0, 1,\infty\}$.
A famous Belyi theorem says
that every smooth projective algebraic curve $X$ defined over the algebraic closure $\overline{\Q}$ (in $\C$) of rational numbers  can be
realized as such a ramified covering of $\CP^1$. Moreover, the universal Galois group $\Gal(\overline{\Q}:\Q)$ acts faithfully on (equivalence classes of)  Belyi maps.
Given a Belyi map $f: X\rar \CP^1$, the associated ribbon hypergraph ({\em dessin d'enfants}) $\Ga$ is embedded into the Riemann surface $X$ as the pre-image
$$
f^{-1}(\xy
 (0,0)*{\bullet}="a",
(5,0)*{\circ}="b",
\ar @{-} "a";"b" <0pt>
\endxy)
$$
of the unit interval $[0,1]\subset \CP^1$ with the point $0$ presented as the white vertex vertex and the point $1$ as the asterisk. For example, the hypergraph
$$
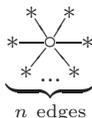
\underbrace{\xy
 (0,0)*{\circ}="a",
(5,0)*{\ast}="b",
(-5,0)*{\ast}="c",
(-2.5,4)*{\ast}="d",
(2.5,4)*{\ast}="e",
(-3.2,-4)*{\ast}="f",
(3.2,-4)*{\ast}="g",
(0,-5)*{...},
\ar @{-} "a";"b" <0pt>
\ar @{-} "a";"c" <0pt>
\ar @{-} "a";"d" <0pt>
\ar @{-} "a";"e" <0pt>
\ar @{-} "a";"f" <0pt>
\ar @{-} "a";"g" <0pt>
\endxy}_{n\ \mathrm{edges}}
$$
corresponds to the Belyi map $f: \CP^1\rar \CP^1$ given by $f(z)=z^n$,
while the hypergraph
$$
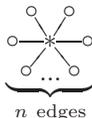
\underbrace{\xy
 (0,0)*{\ast}="a",
(5,0)*{\circ}="b",
(-5,0)*{\circ}="c",
(-2.5,4)*{\circ}="d",
(2.5,4)*{\circ}="e",
(-3.2,-4)*{\circ}="f",
(3.2,-4)*{\circ}="g",
(0,-5)*{...},
\ar @{-} "a";"b" <0pt>
\ar @{-} "a";"c" <0pt>
\ar @{-} "a";"d" <0pt>
\ar @{-} "a";"e" <0pt>
\ar @{-} "a";"f" <0pt>
\ar @{-} "a";"g" <0pt>
\endxy}_{n\ \mathrm{edges}}
$$
corresponds to the Belyi map $f(z)=(1-z)^n$.
However, ribbon hypergraphs
 used in this paper are $\Z$-graded and {\em oriented}\, (see the next subsection) and it is not clear how this extra structure may fit this particular geometric interpretation of hypergraphs.
{\em Oriented}\, ribbon graphs have a useful interpretation as combinatorial objects parameterizing  cells of a certain cell decomposition of the moduli spaces $\cM_{g,n}$ of algebraic curves of genus $g$ with $n$ punctures; they also play an important role in the deformation theory of the Goldman-Turaev Lie bialgebra. We show in this paper that {\em oriented}\, ribbon hypergraphs have much to do with strongly homotopy involutive Lie bialgebras and hence
might be useful in the string topology.

\subsection{Orientation and $\Z$ grading} Let $d$ be an arbitrary integer.
To a ribbon hypergraph $\Ga$ we assign its {\em homological degree}
 $$
 |\Ga|= (d+1)\# H(\Ga) - d\# E(\Ga),
 $$
 i.e. each hyperedge has degree $d+1$ and every edge has degree $-d$.

\sip

 An {\em orientation}
on a  ribbon hypergraph $\Ga$ is, by definition, a
\Bi
\item choice of the total ordering (up to an even permutation) on the set
of hyperedges $H(\Ga)$ {\em for $d$ even},
\item  choice of the compatible (with the underlying cyclic ordering) total ordering --- up to an even permutation --- on each subset of edges $p_\ast^{-1}(h)\subset E(\Ga)$, $\forall h\in H(\Ga)$,  {\em for $d$ odd}.
\Ei
As $E(\Ga)=\sqcup_{h\in H(\Ga)} p_*^{-1}(h)$ and each set  $p_*^{-1}(h)$
is cyclically ordered, a choice of the total ordering in $E(\Ga)$ {\em for $d$ odd} is equivalent to a choice of a total ordering of $H_{odd}(\Ga):=\{h\in H(\Ga)|\ |h|\in 2\Z+1\}$ (up to an even permutation), and a choice of a total ordering of each set $p_*^{-1}(h)$, $h\in H_{even}(\Ga):=\{h\in H(\Ga)|\ |h|\in 2\Z\}$ which is compatible with the given cyclic ordering (again up to an even permutation).

Note that every ribbon hypergraph has precisely two possible orientations. If $\Ga$
is an oriented ribbon hypergraph, then the same hypergraph equipped with an opposite
orientation is denoted by $\Ga^{opp}$.

\sip

Also note that if $\Ga$ has all  hyperedges bivalent, then the above definition
agrees with the notion of orientation in the prop of ribbon graphs $\cR\cG ra_{d+1}$ introduced in [MW].

\subsection{Boundaries and corners of a hypergraph} Let us represent pictorially a ribbon hypergraph $\Ga$ with vertices and hyperedges blown up into dashed and, respectively, double solid circles, for example

$$
 
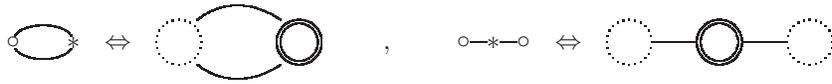
\xy
(-4,0)*{\circ}="1";
(4,0)*{\ast}="2";
"1";"2" **\crv{(-4,3) & (4,3)};
"1";"2" **\crv{(-4,-3) & (4,-3)};
\endxy
\ \ \Leftrightarrow \ \
\xy
 (0,0)*{
\xycircle(3,3){.}};
(16,0)*{
\xycircle(2.5,2.5){-}};
(16,0)*{
\xycircle(3,3){-}};
 (2.5,3)*{}="a",
(13.5,3)*{}="b",
 (2.5,-3)*{}="a'",
(13.5,-3)*{}="b'",
"a";"b" **\crv{(5,5) & (9,6)};
"a'";"b'" **\crv{(5,-5) & (9,-6)};
\endxy\ \ \ \ \ \ \ , \ \ \ \ \ \ \
\xy
 (0,0)*{\ast}="a",
(4,0)*{\circ}="b",
(-4,0)*{\circ}="c",
\ar @{-} "a";"b" <0pt>
\ar @{-} "a";"c" <0pt>
\endxy
\ \ \Leftrightarrow \ \
\xy
 (0,0)*{
\xycircle(3,3){.}};
(12,0)*{
\xycircle(3,3){-}};
(12,0)*{
\xycircle(2.5,2.5){-}};
(24,0)*{
\xycircle(3,3){.}};
 (3,0)*{}="a",
(9,0)*{}="b",
 (21,0)*{}="a'",
 (15,0)*{}="b'",
\ar @{-} "a";"b" <0pt>
\ar @{-} "b'";"a'" <0pt>
\endxy
$$
Edges attached to dashed (resp.\ double solid) circles divide the latters into the disjoint union of chords which can be called {\em vertex (resp.\ hyperedge) corners};
 thus with any vertex $v\in V(\Ga)$ we associate a cyclically ordered set $C(v)$ of its corners (which is, of course,  isomorphic to its set $p_\circ^{-1}(v)$ of attached edges but has a different geometric incarnation).
The motivation for this terminology is that any boundary of a hypergraph can be understood as a polytope glued from edges at that corners. For example, the unique boundary $b$ of the right graph just above is given the following polytope
$$
b=
 \Ba{c}\resizebox{13mm}{!}
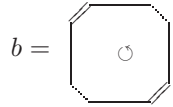
{\xy
(6,-8)*{\circlearrowleft},
 (0,0)*{}="a1",
(10,0)*{}="a2",
(13,-3)*{}="a3",
(13,-13)*{}="a4",
(10,-16)*{}="a5",
(0,-16)*{}="a6",
(-3,-13)*{}="a7",
(-3,-3)*{}="a8",
\ar @{-} "a1";"a2" <0pt>
\ar @{.} "a2";"a3" <0pt>
\ar @{-} "a3";"a4" <0pt>
\ar @{=} "a4";"a5" <0pt>
\ar @{-} "a5";"a6" <0pt>
\ar @{.} "a6";"a7" <0pt>
\ar @{-} "a7";"a8" <0pt>
\ar @{=} "a8";"a1" <0pt>
\endxy}\vspace{2mm}\Ea
$$
where small dashed (resp.\ double solid) intervals stand for the vertex (resp.\ hyperedge) corners. Thus with any boundary $b$ of a hypergraph one can associate a cyclically ordered set $C(b)$ of its vertex corners. We shall use these cyclically ordered sets $C(v)$ and $C(b)$ in the definition of the prop composition of hypergraphs below, and in the construction of canonical representations of that prop in spaces
of cyclic words.

\subsection{Prop of ribbon hypergraphs}  Let $\cR \caH_{m,n}^{k,l}$ be the set of
(isomorphism classes of)
oriented
ribbon hypergraphs $\Ga$ with $n$  vertices labelled by elements of $[n]$,
$k$ edges, $l$ unlabelled hyperedges  and $m$ boundaries labelled by elements of $[m]$. Consider a collection
of quotient $\bS$-bimodules,
$$
\RH_d:=\left\{\RH_d(m,n):=\bigoplus_{k,l\geq 1}
\frac{\K\langle \cR \caH_{m,n}^{k,l} \rangle}{\{ \Ga=-\Ga^{opp}, \Ga\in  \caH_{m,n}^{k,l}\}}
[dk - (d+1)l]\right\}_{m,n\geq 1}
$$
 Thus elements of $\RH$ are isomorphisms classes
of $\Z$-graded oriented ribbon hypergraphs $\Ga$ whose
vertices and boundaries are enumerated. Ribbon hypergraphs admitting automorphisms
which reverse their orientations are equal to zero in $\RH_d$. The $\bS$-module $\RH_d$ contains a submodule $\cR \cG ra_d$ generated by hypergraphs with all hyperedges bivalent. This submodule has a prop structure \cite{MW} which can be easily extended to $\RH_d$. Indeed, the
horizontal composition
$$
 \Ba{rccc}
 \circ: & \RH_d(m_1,n_1)\ot_\K \RH_d(m_2,n_2) &\lon &
 \RH_d(m_1+m_2,n_1+n_2)\\
        &    \Ga_2\ot  \Ga_1 & \lon & \Ga_2\sqcup \Ga_1
        \Ea
 $$
is defined as the disjoint union of ribbon hypergraphs,
and the vertical composition,
 $$
 \Ba{rccc}
 \circ: & \RH_d(p,m)\ot_\K \RH_d(m,n) &\lon & \RH_d(p,n)\\
        &    \Ga_2\ot  \Ga_1 & \lon & \Ga_2\circ \Ga_1
        \Ea
 $$
is defined by gluing, for every  $i\in [m]$, the  $i$-th oriented boundary
$b$ of $\Ga_1$,
$$
b\sim \Ba{c}\resizebox{15mm}{!}
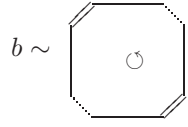
{\xy
(6,-8)*{\circlearrowleft},
 (0,0)*{}="a1",
(10,0)*{}="a2",
(13,-3)*{}="a3",
(13,-13)*{}="a4",
(10,-16)*{}="a5",
(0,-16)*{}="a6",
(-3,-13)*{}="a7",
(-3,-3)*{}="a8",
\ar @{-} "a1";"a2" <0pt>
\ar @{.} "a2";"a3" <0pt>
\ar @{-} "a3";"a4" <0pt>
\ar @{=} "a4";"a5" <0pt>
\ar @{-} "a5";"a6" <0pt>
\ar @{.} "a6";"a7" <0pt>
\ar @{-} "a7";"a8" <0pt>
\ar @{=} "a8";"a1" <0pt>
\endxy}\Ea
$$
with the $i$-th vertex $v$ of $\Ga_2$,
$$
v\sim \Ba{c}\resizebox{19mm}{!}
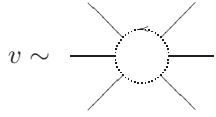
{\xy
 (0,0)*{
\xycircle(6,6){.}};
(-6,0)*{}="1";
(-16,0)*{}="1'";
(4.2,4.2)*{}="2";
(12,12)*{}="2'";
(4.2,-4.2)*{}="3";
(12,-12)*{}="3'";
(6,0)*{}="4";
(16,0)*{}="4'";
(-4.2,4.2)*{}="5";
(-12,12)*{}="5'";
(-4.2,-4.2)*{}="6";
(-12,-12)*{}="6'";
(-0.3,6)*{}="a";
(-0.4,6)*{}="b";
\ar @{-} "1";"1'" <0pt>
\ar @{-} "2";"2'" <0pt>
\ar @{-} "3";"3'" <0pt>
\ar @{-} "4";"4'" <0pt>
\ar @{-} "5";"5'" <0pt>
\ar @{->} "a";"b" <0pt>
\ar @{-} "6";"6'" <0pt>
\endxy}\Ea
$$
by erasing the vertex $v$ from $\Ga_2$ and taking the sum over all possible ways of attaching ``hanging in the air" (half)edges from the set $p_\circ^{-1}(v)$ to the set of  dashed corners from  $C(b)$ while respecting the cyclic structures of both sets; put another way we take a sum over all morphisms $p_\circ^{-1}(v) \rar C(b)$ of cyclically ordered sets.
Every ribbon graph in this linear combination comes equipped naturally with an induced
orientation, and belongs to $\RH_{d}(p,n)$.
The graph $\circ$ consisting of a single white vertex acts as the unit in $\RH_d$.
The subspace of $\RH_{d}$ spanned by {\em connected}\, ribbon graphs forms a
{\em properad}\, which we denote by the same symbol $\RH_{d}$.

\subsubsection{{\bf Theorem}}\label{3: Prop on repr of Rgra in CycW}
{\em Let $W$ be an arbitrary graded vector space and  $
Cyc(W)=\oplus_{n\geq 0} (W^{\ot n})^{\Z_n}
$ the associated space of cyclic words. Then any collection
$$
\Theta_n: (W[d])^{\ot n}_{\Z_n} \lon \K[1+d], \ \ \ \ \ n\geq 1,
$$
of cyclically (skew)invariant maps
gives canonically rise to a representation
$$
\rho_{\Theta_\bu}: \RH_{d} \lon \cE nd_{Cyc(W)}.
$$
of the prop of hypergraphs in  $Cyc(W)$.
}

\begin{proof} If only $\Theta_2$ is non-zero, the associated representation
$$
\rho_{\Theta_2}: \cR \cG ra_d \lon  \cE nd_{Cyc(W)}
$$
was constructed in Theorem~4.2.2 of \cite{MW}. The general case is a straightforward hypergraph extension of that construction. Let us sketch this extension for $d$ even (the case $d$ odd is completely analogous).
Consider a hypergraph $\Ga\in \RH_{d}(m,n)$
 with $n$ vertices $(v_1,\ldots, v_i,\ldots, v_n)$ and $m$ boundaries
 $(b_1,\ldots, b_j,\ldots, b_m)$. Using $\Theta_\bu$ we will construct a linear map
 $$
 \Ba{rccc}
 \rho_{\Theta_\bu}^\Ga: & \ot^n Cyc(W) & \lon & \ot^m Cyc(W)\\
  & \cW_1\ot \ldots \ot \cW_n &\lon &  \rho_{\Theta_\bu}^\Ga( \cW_1, \ldots, \cW_n)
  \Ea
 $$
where
\Beq\label{4: cyclic words cW}
\left\{\cW_i:= (w_{i_1}\ot\ldots\ot w_{{p_i}})_{\Z_{p_i}}\right\}_{ 1\leq i\leq n,\ \ p_i\in \N}\,
\Eeq
is an arbitrary collection of $n$ cyclic words from  $Cyc(W)$. If $\# p_\circ^{-1}(v_i)> p_i$
for at least one $i\in [n]$, i.e.\ if number of edges attached to $v_i$ greater than the length of the word $\cW_i$, we set $\rho_{\Theta_\bu}^\Ga( \cW_1, \ldots, \cW_n)=0$. Otherwise it makes sense to consider a {\em state $s$}\, which is by definition a collection of fixed injective morphisms of cyclically ordered sets
$$
s_i: p_\circ^{-1}(v_i) \lon \{w_{i_1}, \ldots, w_{{p_i}}\},\ \ \forall i\in [n],
$$
that is, an assignment of some letter $w_{i_\bu}$ from the word $\cW_i$ to each edge $e_{i_\bu}\in  p_\circ^{-1}(v_i)$ of each vertex $v_i$ in
a way which respects cyclic orderings of both sets. Note that for each state the complement,  $(w_{i_1},\ldots,  w_{i_{l_p}})\setminus \Img s_i$,
splits into a  disjoint (cyclically ordered)  union of {\em totally ordered}\,  subsets,
$\coprod_{c\in C(v)}I_c$,
    parameterized by the set of corners of the vertex $v$. Note also that to each boundary $b_j\in B(\Ga)$ we can associate
a cyclic word
$$\displaystyle
\cW'_j:=\left(\bigotimes_{c\in C(b_j)} I_c\right)_{\Z_{\sum_{c\in C(b_j)}\# I_c}}
$$
where the tensor product is taken along the given cyclic ordering in the set $C(b_j)$,

\sip

Recall that the set of hyperedges $H(\Ga)$ is defined as the set of orbits $E(\Ga)/\sigma_1$ of the permutation
$\sigma_1$. To any hyperedge   $h\in H(\Ga)$ of valency $|h|\in \N$ there corresponds
therefore a cyclically ordered set of $|h|$ edges $p_*^{-1}(h)\subset E(\Ga)$. Let us choose for a moment a compatible total order on this set, i.e.\ write it as
$$
p_{\ast}^{-1}(h)=\{e^h_1, e^h_2=\sigma_1(e^h_1), \ldots, \ldots e^h_{|h|}=\sigma_1^{|h|-1}(e^h_1)\}
$$
 for some chosen edge $e^h_1\in p_\ast^{-1}(h)$. As the set of edges decomposes into
the disjoint union
$$
E(\Ga)=\coprod_{h\in H(\Ga)} p_*^{-1}(h),
$$
we can use the given maps $\Theta_k: (\ot^k W)_{\Z_k} \rar \K[1+d-dk]$ to define the {\em weight}\,
of any given state $s$ on $\Ga$ as the following number,
    $$
    \la_s:=\prod_{h\in H(\Ga)}\Theta_{|h|} \left(s(e^h_1), \ldots s(e^h_{|h|})\right).
    $$

    \sip

Thus to each state $s=\{s_i\}_{i\in [n]}$ we
associate an element
$$
\rho_{\Theta_\bu}^s(\cW_1,\ldots, \cW_n):=
(-1)^\sigma \la_s
\cW'_{b_1}\ot\ldots \ot \cW'_{b_m} \in \ot^m Cyc(W)
$$
where $(-1)^\sigma$ is the standard Koszul sign of the regrouping permutation,
$$
\sigma:\cW_1\ot\ldots \ot \cW_n\lon \prod_{h\in  (\Ga)}
 \left(s(e^h_1), \ldots s(e^h_{|h|})\right)
\ot \cW'_{b_1}\ot\ldots \ot \cW'_{b_m}.
$$
Note that $\rho_{\Theta_\bu}^s(\cW_1,\ldots, \cW_n)$ does not depend on the choices of compatible total orderings in the sets
  $p_\ast^{-1}(h)$ made above.
\sip

Finally we define a linear map,
$$
\Ba{rccc}
\rho_{\Theta_\bu}: & \cR\cG ra_{d}(m,n) & \lon & \Hom(\ot^n Cyc(W), \ot^m Cyc(W)) \\
 & \Ga & \lon & \rho_{\Theta_\bu}^\Ga
 \Ea
$$
by setting the value of $ \rho_{\Theta_\bu}^\Ga$ on cyclic words (\ref{4: cyclic words cW}) to be
equal to
$$
\rho_{\Theta_\bu}^\Ga(\cW_1,\ldots, \cW_n):=\left\{\Ba{cl} 0 & \mbox{if}\ \# p_\circ^{-1}(v_i)> p_i\
\mbox{for some}\ i\in [n] \\
\displaystyle \sum_{\mathrm{all\ possiblle}\atop \mathrm{states}\ s} \rho_{\Theta_\bu}^s(\cW_1,\ldots,
 \cW_n) &  \mbox{otherwise} \Ea\right.
$$
It is now straightforward to check that the map $ \rho_{\Theta_\bu}$ respects prop
compositions  in  $\RH_d$ and $\cE nd_{Cyc(W)}$ because the prop structure in
former has been just read off from the compositions of operators  $\rho_{\Theta_\bu}^\Ga$ in the latter.
\end{proof}

\bip


{\large
\section{\bf Strongly homotopy involutive Lie bialgebras and ribbon hypergraphs}
}

\subsection{{\bf Reminder from \cite{MW}}}
{\em There is a morphism of props,
$$
\rho: \LoBd \lon \RH_d
$$
given on generators as follows,
$$
\rho\left(
 \begin{xy}
 <0mm,-0.55mm>*{};<0mm,-2.5mm>*{}**@{-},
 <0.5mm,0.5mm>*{};<2.2mm,2.2mm>*{}**@{-},
 <-0.48mm,0.48mm>*{};<-2.2mm,2.2mm>*{}**@{-},
 <0mm,0mm>*{\circ};<0mm,0mm>*{}**@{},
 \end{xy}\right)=\xy
(-4,0)*{\circ}="1";
(4,0)*{\ast}="2";
"1";"2" **\crv{(-4,3) & (4,3)};
"1";"2" **\crv{(-4,-3) & (4,-3)};
\endxy
 \ \ \  , \ \ \
\rho\left(
 \begin{xy}
 <0mm,0.66mm>*{};<0mm,3mm>*{}**@{-},
 <0.39mm,-0.39mm>*{};<2.2mm,-2.2mm>*{}**@{-},
 <-0.35mm,-0.35mm>*{};<-2.2mm,-2.2mm>*{}**@{-},
 <0mm,0mm>*{\circ};<0mm,0mm>*{}**@{},
 \end{xy}\right)=\xy
 (0,0)*{\ast}="a",
(4,0)*{\circ}="b",
(-4,0)*{\circ}="c",
\ar @{-} "a";"b" <0pt>
\ar @{-} "a";"c" <0pt>
\endxy
$$
}

\sip

The main new result of this note
 is an observation that this maps lifts to a morphism of dg props $\HoLoB_d\rar \RH_d$ which is non-trivial on all generators (see below).

\subsection{Proposition} {\em
There is a morphism of dg props,
$$
\rho: \HoLBd \lon \RH_d
$$
given on generators as follows,
$$
\Ba{c}\resizebox{14mm}{!}{\begin{xy}
 <0mm,0mm>*{\circ};<0mm,0mm>*{}**@{},
 <-0.6mm,0.44mm>*{};<-8mm,5mm>*{}**@{-},
 <-0.4mm,0.7mm>*{};<-4.5mm,5mm>*{}**@{-},
 <0mm,0mm>*{};<-1mm,5mm>*{\ldots}**@{},
 <0.4mm,0.7mm>*{};<4.5mm,5mm>*{}**@{-},
 <0.6mm,0.44mm>*{};<8mm,5mm>*{}**@{-},
   <0mm,0mm>*{};<-8.5mm,5.5mm>*{^1}**@{},
   <0mm,0mm>*{};<-5mm,5.5mm>*{^2}**@{},
   <0mm,0mm>*{};<4.5mm,5.5mm>*{^{m\hspace{-0.5mm}-\hspace{-0.5mm}1}}**@{},
   <0mm,0mm>*{};<9.0mm,5.5mm>*{^m}**@{},
 <-0.6mm,-0.44mm>*{};<-8mm,-5mm>*{}**@{-},
 <-0.4mm,-0.7mm>*{};<-4.5mm,-5mm>*{}**@{-},
 <0mm,0mm>*{};<-1mm,-5mm>*{\ldots}**@{},
 <0.4mm,-0.7mm>*{};<4.5mm,-5mm>*{}**@{-},
 <0.6mm,-0.44mm>*{};<8mm,-5mm>*{}**@{-},
   <0mm,0mm>*{};<-8.5mm,-6.9mm>*{^1}**@{},
   <0mm,0mm>*{};<-5mm,-6.9mm>*{^2}**@{},
   <0mm,0mm>*{};<4.5mm,-6.9mm>*{^{n\hspace{-0.5mm}-\hspace{-0.5mm}1}}**@{},
   <0mm,0mm>*{};<9.0mm,-6.9mm>*{^n}**@{},
 \end{xy}}\Ea \stackrel{\rho}{\lon} \sum
 \underbrace{\xy
 <5mm,-5mm>*{...};
(0,5)*{\ast}="1";
(-7,-5)*{\circ}="2";
(-1,-5)*{\circ}="3";
(9,-5)*{\circ}="4";
"1";"2" **\crv{(0,4.5) & (0,0)};
"1";"2" **\crv{(-5,0) & (-4,1)};
"1";"3" **\crv{(0,5.5) & (2,4)};
"1";"4" **\crv{(0,4.5) & (2,7)};
\ar @{-} "1";"2" <0pt>
\ar @{-} "1";"3" <0pt>
\ar @{-} "1";"4" <0pt>
\endxy}_{n+m-1\ edges\atop m\ boundaries}
$$
where the sum on the right hand side is over all possible ways of attaching
$n+m-1$ edges beginning at the asterisk vertex to $n$ white vertices
(whose numerical labels are (skew)symmetrized) in such a way that every white vertex is hit and the total number of boundaries of the resulting hypergraph
equals precisely $m$ (and their numerical labels are also (skew)symmetrized).
}

\mip

\begin{proof} As the prop $\RH_d$ has vanishing differential, the Propsoition holds true if and only of
$$
\rho\left(
\delta
\Ba{c}\resizebox{14mm}{!}{\begin{xy}
 <0mm,0mm>*{\circ};<0mm,0mm>*{}**@{},
 <-0.6mm,0.44mm>*{};<-8mm,5mm>*{}**@{-},
 <-0.4mm,0.7mm>*{};<-4.5mm,5mm>*{}**@{-},
 <0mm,0mm>*{};<-1mm,5mm>*{\ldots}**@{},
 <0.4mm,0.7mm>*{};<4.5mm,5mm>*{}**@{-},
 <0.6mm,0.44mm>*{};<8mm,5mm>*{}**@{-},
   <0mm,0mm>*{};<-8.5mm,5.5mm>*{^1}**@{},
   <0mm,0mm>*{};<-5mm,5.5mm>*{^2}**@{},
   <0mm,0mm>*{};<4.5mm,5.5mm>*{^{m\hspace{-0.5mm}-\hspace{-0.5mm}1}}**@{},
   <0mm,0mm>*{};<9.0mm,5.5mm>*{^m}**@{},
 <-0.6mm,-0.44mm>*{};<-8mm,-5mm>*{}**@{-},
 <-0.4mm,-0.7mm>*{};<-4.5mm,-5mm>*{}**@{-},
 <0mm,0mm>*{};<-1mm,-5mm>*{\ldots}**@{},
 <0.4mm,-0.7mm>*{};<4.5mm,-5mm>*{}**@{-},
 <0.6mm,-0.44mm>*{};<8mm,-5mm>*{}**@{-},
   <0mm,0mm>*{};<-8.5mm,-6.9mm>*{^1}**@{},
   <0mm,0mm>*{};<-5mm,-6.9mm>*{^2}**@{},
   <0mm,0mm>*{};<4.5mm,-6.9mm>*{^{n\hspace{-0.5mm}-\hspace{-0.5mm}1}}**@{},
   <0mm,0mm>*{};<9.0mm,-6.9mm>*{^n}**@{},
 \end{xy}}\Ea\right)
\ \ = \ \ \rho
\left(
 \sum_{[1,\ldots,m]=I_1\sqcup I_2\atop
 {|I_1|\geq 0, |I_2|\geq 1}}
 \sum_{[1,\ldots,n]=J_1\sqcup J_2\atop
 {|J_1|\geq 1, |J_2|\geq 1}
}\hspace{0mm}
\pm
\Ba{c}\resizebox{22mm}{!}{ \begin{xy}
 <0mm,0mm>*{\circ};<0mm,0mm>*{}**@{},
 <-0.6mm,0.44mm>*{};<-8mm,5mm>*{}**@{-},
 <-0.4mm,0.7mm>*{};<-4.5mm,5mm>*{}**@{-},
 <0mm,0mm>*{};<0mm,5mm>*{\ldots}**@{},
 <0.4mm,0.7mm>*{};<4.5mm,5mm>*{}**@{-},
 <0.6mm,0.44mm>*{};<12.4mm,4.8mm>*{}**@{-},
     <0mm,0mm>*{};<-2mm,7mm>*{\overbrace{\ \ \ \ \ \ \ \ \ \ \ \ }}**@{},
     <0mm,0mm>*{};<-2mm,9mm>*{^{I_1}}**@{},
 <-0.6mm,-0.44mm>*{};<-8mm,-5mm>*{}**@{-},
 <-0.4mm,-0.7mm>*{};<-4.5mm,-5mm>*{}**@{-},
 <0mm,0mm>*{};<-1mm,-5mm>*{\ldots}**@{},
 <0.4mm,-0.7mm>*{};<4.5mm,-5mm>*{}**@{-},
 <0.6mm,-0.44mm>*{};<8mm,-5mm>*{}**@{-},
      <0mm,0mm>*{};<0mm,-7mm>*{\underbrace{\ \ \ \ \ \ \ \ \ \ \ \ \ \ \
      }}**@{},
      <0mm,0mm>*{};<0mm,-10.6mm>*{_{J_1}}**@{},
 <13mm,5mm>*{};<13mm,5mm>*{\circ}**@{},
 <12.6mm,5.44mm>*{};<5mm,10mm>*{}**@{-},
 <12.6mm,5.7mm>*{};<8.5mm,10mm>*{}**@{-},
 <13mm,5mm>*{};<13mm,10mm>*{\ldots}**@{},
 <13.4mm,5.7mm>*{};<16.5mm,10mm>*{}**@{-},
 <13.6mm,5.44mm>*{};<20mm,10mm>*{}**@{-},
      <13mm,5mm>*{};<13mm,12mm>*{\overbrace{\ \ \ \ \ \ \ \ \ \ \ \ \ \ }}**@{},
      <13mm,5mm>*{};<13mm,14mm>*{^{I_2}}**@{},
 <12.4mm,4.3mm>*{};<8mm,0mm>*{}**@{-},
 <12.6mm,4.3mm>*{};<12mm,0mm>*{\ldots}**@{},
 <13.4mm,4.5mm>*{};<16.5mm,0mm>*{}**@{-},
 <13.6mm,4.8mm>*{};<20mm,0mm>*{}**@{-},
     <13mm,5mm>*{};<14.3mm,-2mm>*{\underbrace{\ \ \ \ \ \ \ \ \ \ \ }}**@{},
     <13mm,5mm>*{};<14.3mm,-4.5mm>*{_{J_2}}**@{},
 \end{xy}}\Ea\right)= 0.
$$
This is almsot obvious as the r.h.s. is given by the sum
$$
\sum_{[1,\ldots,m]=I_1\sqcup I_2\atop
 {|I_1|\geq 0, |I_2|\geq 1}}
 \sum_{[1,\ldots,n]=J_1\sqcup J_2\atop
 {|J_1|\geq 1, |J_2|\geq 1}
}
\sum
\pm
 \underbrace{\xy
 <5mm,-5mm>*{...};
(0,5)*{\ast}="1";
(-7,-5)*{\circ}="2";
(-1,-5)*{\circ}="3";
(9,-5)*{\circ}="4";
<20mm,-5mm>*{...};
(16,5)*{\ast}="1'";
(9,-5)*{\circ}="2'";
(15,-5)*{\circ}="3'";
(24,-5)*{\circ}="4'";
"1";"2" **\crv{(0,4.5) & (0,0)};
"1";"2" **\crv{(-5,0) & (-4,1)};
"1";"3" **\crv{(0,5.5) & (2,4)};
"1";"4" **\crv{(0,4.5) & (2,7)};
"1'";"2'" **\crv{(15,4.5) & (15,0)};
"1'";"2'" **\crv{(11,0) & (12,1)};
"1'";"3'" **\crv{(16,5.5) & (18,4)};
"1'";"4'" **\crv{(16,4.5) & (18,7)};
\ar @{-} "1";"2" <0pt>
\ar @{-} "1";"3" <0pt>
\ar @{-} "1";"4" <0pt>
\ar @{-} "1'";"3'" <0pt>
\ar @{-} "1'";"4'" <0pt>
\endxy}_{\ \ \ \ \ \ \ |J_1|+|I_1|\ \ +  \ \ |J_2|+|I_2| edges \atop \ \ \ \ \ \ \ \ \ |I_1| \ \ \ \ \ \ + \ \ \ \ |I_2| boundaries}
$$
 which
 vanishes in $\RH_d$ for symmetry reasons (it is most easy to check this claim in the case $d$ is even when the asterisk vertices are odd, and the symbol $\pm$ becomes $+$).
\end{proof}

\mip

There is a canonical morphism of props $\HoLBd\rar \HoLoBd$. The above morphism factors through the composition
$\HoLBd\rar \HoLoBd \stackrel{\rho^\diamond}{\lon} \RH_d$.

\subsection{Theorem}\label{4: Main Th}  {\em
There is a morphism of dg props,
$$
\rho^\diamond: \HoLoBd \lon \RH_d
$$
given on generators as follows,
\Beq\label{4: rho-diamond}
\Ba{c}\resizebox{16mm}{!}{\xy
(-9,-6)*{};
(0,0)*+{a}*\cir{}
**\dir{-};
(-5,-6)*{};
(0,0)*+{a}*\cir{}
**\dir{-};
(9,-6)*{};
(0,0)*+{a}*\cir{}
**\dir{-};
(5,-6)*{};
(0,0)*+{a}*\cir{}
**\dir{-};
(0,-6)*{\ldots};
(-10,-8)*{_1};
(-6,-8)*{_2};
(10,-8)*{_n};
(-9,6)*{};
(0,0)*+{a}*\cir{}
**\dir{-};
(-5,6)*{};
(0,0)*+{a}*\cir{}
**\dir{-};
(9,6)*{};
(0,0)*+{a}*\cir{}
**\dir{-};
(5,6)*{};
(0,0)*+{a}*\cir{}
**\dir{-};
(0,6)*{\ldots};
(-10,8)*{_1};
(-6,8)*{_2};
(10,8)*{_m};
\endxy}\Ea
\stackrel{\rho^\diamond}{\lon}
 \sum
 \underbrace{\xy
 <2mm,-5mm>*{...};
(0,5)*{\ast}="1";
(-7,-5)*{\circ}="2";
(-1,-5)*{\circ}="3";
(9,-5)*{\circ}="4";
"1";"2" **\crv{(0,4.5) & (0,0)};
"1";"2" **\crv{(-5,0) & (-4,1)};
"1";"2" **\crv{(-14,-4) & (0,0)};
"1";"3" **\crv{(0,5.5) & (2,4)};
"1";"4" **\crv{(10,4.5) & (-2,-7)};
\ar @{-} "1";"3" <0pt>
\ar @{-} "1";"4" <0pt>
\endxy}_{n+m+2a-1\ edges\atop m\ boundaries}
\Eeq
where the sum on the right hand side is over all possible ways to attach
$n+m+2a-1$ edges beginning at the asterisk vertex to $n$ white vertices
(whose numerical labels are (skew)symmetrized) in such a way that every white vertex is hit and the total number of boundaries of the resulting hypergraph
equals precisely $m$  (and their numerical labels are (skew)symmetrized).

}

\sip

The orientations of the hypergraphs shown in the above formula are determined uniquely by a simple $\caH o\cB\cV^{com}_d$ operator which is constructed in the proof.

\begin{proof}
The composition of the above map $\rho^\diamond$ with the canonical representation
$\rho_{\Theta_\bu}$ from Proposition {\ref{3: Prop on repr of Rgra in CycW}} implies that for any collection of cyclically (skew)symmetric maps \ref{1: Theta_n} on a graded vector space $W$  the associated  vector space $Cyc(W)$ is canonically a $\HoLoBd$ algebra. In fact one can read the Theorem from this conclusion provided the latter is independent of choices of $W$ and the higher products $\Theta_\bu$.

\bip

  Assume $d$ even. To make the construction of the representation $\rho_{\Theta_\bu}\circ \rho^\diamond$ as simple and transparent as possible  we shall employ, following Serguei Barannikov \cite{Ba3}, the  invariant theory and identify the space
$Cyc(W)$ with the space of $GL(\K^N)$-invariants,
$$
Cyc(W)=\lim_{N\rar \infty} \oplus_{n\geq 0} \left(\ot^n (W\ot \End(\K^n)\right)^{GL(\K^N)}
,
$$
that is, we interpret a cyclic word, $\cW=(w_{a_1}\ot\ldots\ot w_{a_n})^{\Z_n}$,
$\{w_i\}_{i\in I}$ a basis in $W$, with the trace of the product of $N\times N$ matrices,
$$
\cW=tr\left(A_{a_1}A_{a_2}\cdots A_{a_n} \right)=\sum_{\al_\bu} A_{a_1\, \al_1}^{\al_0} A_{a_2\, \al_2}^{\al_1}\cdots A_{a_n\, \al_0}^{\al_n}
, \ \ \ \ \ A_{i_a}:=\left(A_{i_a\, \be}^{\al}\right) \in \End(\K^N), \ \ a\in [n],\ \al,\be\in [N],
$$
for sufficiently large $N\in \N$. The graded cyclically symmetric maps
$$
\Ba{rccc}
\Theta_n: & \ot^n W & \lon &\K[1+d-nd] \\
        & w_{a_1}\ot ...\ot  w_{a_n} & \lon & \Theta_{a_1...a_n}
        \Ea
$$
define a degree $1$  operator
$$
\Delta=\sum_{n\geq 1} \hbar^{n-1}  \Theta_{a_1...a_n}\frac{\p^n}{\p A_{a_1\, \al_1}^{\ \al_0} \p A_{a_2\, \al_2}^{\ \al_1}\cdots \p A_{a_n\, \al_0}^{\ \al_n}}
$$
on $\odot^\bu (Cyc(W)[-d])[[\hbar]]$ whose square is obviously zero; here the formal parameter $\hbar$ has degree $2d$. The latter defines
a $\caH o\cB\cV_d^{com}$-structure in $\odot^\bu (Cyc(W)[-d])[\hbar]]$ and, as explained in \S {\ref{2: Com BV to Holieb}},  an associated $\HoLoBd\{d\}$ algebra structure
in $Cyc(W)[-d]$ which in turn defines a $\HoLoBd$-structure in $Cyc(W)$.

\sip

If $d$ is odd, one can again use a trick from \cite{Ba3} which replaces the ordinary trace of the standard matrix superalgebra with the {\em odd trace}\, of the Bernstein-Leites matrix sub-superalgebra. In fact the construction in \cite{Ba3} explains the construction of the $\caH o\cB\cV_d^{com}$ operator in the case when only $\Theta_2$ is non-zero, and its extension to the general case (\ref{1: Theta_n}) is completely analogous to the $d$ even case discussed  above.
\end{proof}

\subsection{Corollary}\label{4: Corollary about HoLoBd}
{\em Let $W$ be a graded vector space equipped with a collection of linear maps
\Beq\label{4: Theta_n}
\Ba{rccc}
\Theta_p: & \ot^p(W[d])_{\Z_p} & \lon & \K[1+d], \ \ \ \ p\geq 2,
\Ea
\Eeq
Each map $\Theta_p$ makes the (reduced) graded vector space of cyclic words in $W$,
$$
Cyc(W)=\oplus_{n\geq 0} (\ot^n W)_{\Z_n}\ \  \text{or}\ \  \overline{Cyc}(W)=\oplus_{n\geq 1} (\ot^n W)_{\Z_n},
$$
into a
$\HoLoBd$-algebra with only those $\HoLoBd$-operations
$$
\mu(m,n,a)\simeq
\Ba{c}\resizebox{14mm}{!}{\begin{xy}
 <0mm,0mm>*{\circ};<0mm,0mm>*{}**@{},
 <-0.6mm,0.44mm>*{};<-8mm,5mm>*{}**@{-},
 <-0.4mm,0.7mm>*{};<-4.5mm,5mm>*{}**@{-},
 <0mm,0mm>*{};<-1mm,5mm>*{\ldots}**@{},
 <0.4mm,0.7mm>*{};<4.5mm,5mm>*{}**@{-},
 <0.6mm,0.44mm>*{};<8mm,5mm>*{}**@{-},
   <0mm,0mm>*{};<-8.5mm,5.5mm>*{^1}**@{},
   <0mm,0mm>*{};<-5mm,5.5mm>*{^2}**@{},
   <0mm,0mm>*{};<4.5mm,5.5mm>*{^{m\hspace{-0.5mm}-\hspace{-0.5mm}1}}**@{},
   <0mm,0mm>*{};<9.0mm,5.5mm>*{^m}**@{},
 <-0.6mm,-0.44mm>*{};<-8mm,-5mm>*{}**@{-},
 <-0.4mm,-0.7mm>*{};<-4.5mm,-5mm>*{}**@{-},
 <0mm,0mm>*{};<-1mm,-5mm>*{\ldots}**@{},
 <0.4mm,-0.7mm>*{};<4.5mm,-5mm>*{}**@{-},
 <0.6mm,-0.44mm>*{};<8mm,-5mm>*{}**@{-},
   <0mm,0mm>*{};<-8.5mm,-6.9mm>*{^1}**@{},
   <0mm,0mm>*{};<-5mm,-6.9mm>*{^2}**@{},
   <0mm,0mm>*{};<4.5mm,-6.9mm>*{^{n\hspace{-0.5mm}-\hspace{-0.5mm}1}}**@{},
   <0mm,0mm>*{};<9.0mm,-6.9mm>*{^n}**@{},
 \end{xy}}\Ea: \ot^n Cyc(W) \lon \ot^m Cyc(W)
 , \ \ \ m\geq1,n\geq 1, m+n+2a\geq 3,
 $$
 non-trivial which satisfy the condition $m+n+2a=p+1$.
 Moreover, this $\HoLoBd$-algebra structure admits an independent rescaling automorphism
$$
\mu(m,n, a) \rar \la_{m+n+2a}\cdot  \mu(n,m,a), \ \ \ \la_{m+n+2a}\in \K^*,
$$
for each fixes sum  $m+n+2a$ of the integer parameters.

\sip

In the case $p=2$ this construction reproduces  the classical construction of $\LoBd$-algebra structure on $Cyc(W)$ given in \S {\ref{2: subsection on cyclic words}}.}

\subsection{Rescaling freedom} Note that each map $\Theta_k$ from the family
(\ref{1: Theta_n}) can be independently rescaled, $\Theta_k \rar \la_k \Theta_k$, $\forall \la_n\in \K$, so that the morphism $\rho^\diamond$ in (\ref{4: rho-diamond}) can also be rescaled by infinitely many independent parameters --- just rescale in that formula each hyperedge $h$ of valency $k$  by
$$
 h=\xy
 (0,0)*{\ast}="a",
(4,0)*{}="b1",
(3,-3)*{}="b2",
(0,-4)*{}="b3",
(-3,-3)*{}="b4",
(-4,0)*{}="b5",
(-3,3)*{}="b6",
(0,4)*{}="b7",
(3,3)*{}="b8",
\ar @{-} "a";"b1" <0pt>
\ar @{-} "a";"b2" <0pt>
\ar @{-} "a";"b3" <0pt>
\ar @{-} "a";"b4" <0pt>
\ar @{-} "a";"b5" <0pt>
\ar @{-} "a";"b6" <0pt>
\ar @{-} "a";"b7" <0pt>
\ar @{-} "a";"b8" <0pt>
\endxy \ \
\lon \ \
\la_k \cdot \
 \xy
 (0,0)*{\ast}="a",
(4,0)*{}="b1",
(3,-3)*{}="b2",
(0,-4)*{}="b3",
(-3,-3)*{}="b4",
(-4,0)*{}="b5",
(-3,3)*{}="b6",
(0,4)*{}="b7",
(3,3)*{}="b8",
\ar @{-} "a";"b1" <0pt>
\ar @{-} "a";"b2" <0pt>
\ar @{-} "a";"b3" <0pt>
\ar @{-} "a";"b4" <0pt>
\ar @{-} "a";"b5" <0pt>
\ar @{-} "a";"b6" <0pt>
\ar @{-} "a";"b7" <0pt>
\ar @{-} "a";"b8" <0pt>
\endxy\ \ \ \ ,\ \ \ \forall \la_k\in \K,
$$
and get a new morphism $\rho_{\la_\bu}^\diamond$ from $\HoLoBd$ to $\RH_d$.
Such a phenomenon occurs in the string topology --- see Theorem 6.2 and Corollary 6.3 in \cite{CS} --- and its main technical origin in our case is that the prop $\RH_d$ has vanishing differential.

\bip


{\large
\section{\bf An algebraic application: a new family of combinatorial $\HoLoB$-algebras}
}

Given any collection of formal letters $e_1,\ldots, e_N$, i.e.\ given any natural number $N\geq 1$, there is an associated involutive
Lie bialgebra structure on the vector space $Cyc(W_N)$ of cyclic words,
$$
W_N:= \text{span}_\K \langle e_1,e_2, \ldots, e_N \rangle
$$
which belongs to the family of combinatorial $\LoB$-algebras constructed by Schedler in \cite{Sch} out of any quiver. This particular $\LoB$-algebra has an important geometric meaning (discussed in  \S 1) and corresponds to the quiver
(\ref{1: quiver for W_N}).

\sip

In this section we use Theorem~{\ref{4: Main Th}} and Corollary {\ref{4: Corollary about HoLoBd}} in the case $d=1$ to extend that particular Schedler's construction
to a highly non-trivial (i.e.\ with {\em all}\, higher homotopy operations non-zero) $\HoLoB$-algebra structure on the vector space $Cyc(\widehat{W}_N)$ of cyclic words generated by {\em $\Z$-graded}\, formal letters,
$$
\widehat{W}_N:= \text{span}_\K \langle e_1[-p], \ldots, e_N[-p] \rangle_{p\in \N}
=  \bigoplus_{p\geq 0} W_N[-p]
$$
 where $e_\al[-p]$, $\al\in [N]$, stands for the copy of the formal letter $e_\al$ to which we assign the homological degree $p$.
  Note that $\widehat{W}_N$ has no natural higher products (\ref{1: Theta_n}) so we can not apply Corollary~{\ref{4: Corollary about HoLoBd}} immediately.
 The idea of our construction is to inject first
$$
\Ba{rccc}
u: &\widehat{W}_N & \lon & \widehat{\bW}_N\\
  \Ea
$$
$\widehat{W}_N$ into a larger space $\widehat{\bW}_N$ which does have a natural family  of cyclically (skew)symmetric higher products
$$
\left\{\Theta_{k+2}: \left(\ot^{k+2} (\widehat{\bW}_N[1])\right)_{\Z_{k+2}} \lon \K[2]    \right\}_{k\geq 0}
$$
 so that the associated vector space of cyclic words $Cyc(\widehat{\bW}_N)$  comes equipped with a canonical $\HoLoB$-algebra structure by Corollary {\ref{4: Corollary about HoLoBd}}. The second step will be to check that the image of $u$ is closed with respect to all strongly homotopy involutive Lie bialgebra operations. To realize this programme, consider, for any integer $p\geq 0$,  a set of  $p+2$ copies of the vector space $W_N$
$$
W_N^{0_p}:=W_N[-p]\ ,\ W_N^{1_p}:=W_N\ , \ldots,\  W_N^{{p+1}_p}:=W_N\ ,
$$
one of them (say, labelled by zero)  assigned a shifted homological degree, and define
$$
\widehat{\bW}_N:=\bigoplus_{p\geq 0} W_{N,p}, \ \ \ \ \  W_{N,p}:=W_N^{0_p}\oplus W_N^{1_p}\oplus \ldots \oplus \ldots W_N^{{p+1}_p}
$$
 The vector space $\widehat{\bW}_N$
is countably dimensional, and is equipped by construction with a distinguished basis
$$
\left\{ e^{l_p}_\al\right\}_{p\geq 0, 0\leq  l\leq {p+1}, 1\leq \al \leq N}
$$
where
 $\{ e^{l_p}_\al\}_{\al\in [N]}$ stands for the standard basis of the copy $W^{l_p}_N$. Note that for any $\al\in N$ the basis vector  $e^{l_p}_\al$ has homological degree $0$ if $l\geq 1$, and $-p$ of $l=0$.
  Note also that summands  $W_N^{l_{p'}}$ and
$W_N^{l_{p''}}$  in $\widehat{\bW}_N$ are viewed as different copies of $W_N$ for $p'\neq p''$
(as they belong to different vector spaces $W_{N,p'}$ and $W_{N,p''}$).

\sip

Let us introduce next an infinite family
of cyclically
 (skew)symmetric higher products (the case $d=1$ in the notation (\ref{1: Theta_n}))
$$
\Theta_{k+2}: \ot^{ k+2 }\hat{\bW}_N \lon \K[-k], \ \ \ \ \ \forall \ k\geq 0,
$$
by setting
$$
\Theta_{k+2}\left(
e^{l^1_{\ p_1}}_{\al_1}\ot \ldots \ot  e^{l^{k+2}_{\ \ \ p_{k+2}}}_{\al_{k+2}}
\right):=0
$$
unless
\Bi
\item[(i)] $\al_1=\al_2=\ldots=\al_{k+1}$;
\item[(ii)] $p_1=p_2=\ldots=p_{k+2}=k$;
\item[(iii)]

$
(l^1, \ldots, l^{k+2})=(j, j-1, \ldots, 2,1,0, k+1,k,\ldots,j+1)$  for some $j\in \{0,1,2,\ldots, k+1\}$, i.e.\ we have an isomorphism of cyclically  ordered sets
 $$
 \Ba{c}\resizebox{28mm}{!}
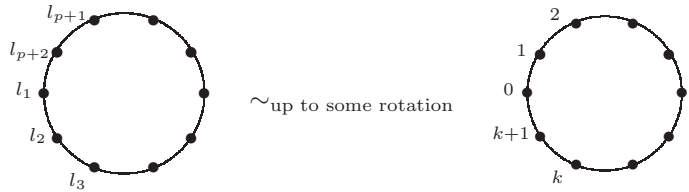
{
\begin{xy}
*\xycircle<30pt,30pt>{};
 <-30pt,0mm>*{\bu};
  <-25pt,15pt>*{\bullet};
<-25pt,15pt>*{\bullet};
<30pt,0mm>*{\bullet};
  <-25pt,-17pt>*{\bullet};
<25pt,-17pt>*{\bullet};
  <25pt,15pt>*{\bullet};
<25pt,15pt>*{\bullet};
  <-11pt,27pt>*{\bullet};
<11pt,27pt>*{\bullet};
  <-11pt,-28pt>*{\bullet};
<11pt,-28pt>*{\bullet};
<-37pt,0pt>*{^{l_{1}}};
<-35pt,15pt>*{^{l_{p+2}}};
<-32pt,-17pt>*{^{l_2}};
<-21pt,29pt>*{^{l_{p+1}}};
<-17pt,-34pt>*{^{l_{3}}};
\end{xy}}\Ea \ \  \sim_{\text{up to some rotation}} \ \
\Ba{c}\resizebox{26mm}{!}{\begin{xy}
*\xycircle<30pt,30pt>{};
 <-30pt,0mm>*{\bu};
  <-25pt,15pt>*{\bullet};
<-25pt,15pt>*{\bullet};
<30pt,0mm>*{\bullet};
  <-25pt,-17pt>*{\bullet};
<25pt,-17pt>*{\bullet};
  <25pt,15pt>*{\bullet};
<25pt,15pt>*{\bullet};
  <-11pt,27pt>*{\bullet};
<11pt,27pt>*{\bullet};
  <-11pt,-28pt>*{\bullet};
<11pt,-28pt>*{\bullet};
<-37pt,0pt>*{^{0}};
<-32pt,15pt>*{^{1}};
<-36pt,-17pt>*{^{k+1}};
<-19pt,29pt>*{^{2}};
<-18pt,-34pt>*{^{{k}}};
\end{xy}}\Ea
$$
\Ei
If all the above conditions are satisfied, we set
$$
\Theta_{k+2}\left(e_\al^{j_k}\ot
e_\al^{j-1_k}\ot \dots e_\al^{1_k}\ot e_\al^{0_p}\ot e_\al^{k+1_k} \ot
e_\al^{{k}_k}\ot \ldots \ot e_\al^{j+1_k}\right):=(-1)^{j(k+1)}
$$
By Corollary {\ref{4: Corollary about HoLoBd}}, the graded vector space $Cyc^\bu(\widehat{\bW}_N)$ is a $\HoLoB$-algebra equipped with quite explicit strongly homotopy operations.
There is a canonical (homogeneous of homological degree zero) injection
$$
\Ba{rccc}
u: & Cyc^\bu(\bW_N) & \lon & Cyc^\bu (\widehat{\bW}_N) \\
& \left( e_{\al_1}[-p_1]\ot \ldots \ot e_{\al_n}[-p_n]\right)^{\Z_n}
& \lon &  u\left(  e_{\al_1}[-p_1]\ot \ldots \ot e_{\al_n}[-p_n]\right)^{\Z_n}
\Ea
$$
identifying each letter $e_\al[-p]$ in a cyclic word from $Cyc(\bW_N)$
with a (totally ordered) word in $p+2$ letters
$$
e_\al[-p] \lon e_\al^{(0_p)} \ot e_\al^{(1_p)}\ot \ldots \ot e_\al^{({p+1}_p)}
$$
i.e.
$$
u\left(  e_{\al_1}[-p_1]\ot ... \ot e_{\al_n}[-p_n]\right)^{\Z_n}:=
\left(e_{\al_1}^{(0)_{p_1}} \ot e_{\al_1}^{(1)_{p_1}}\ot ... \ot e_{\al_1}^{(p_1+1)_{p_1}}\ot \ldots \ot e_{\al_n}^{(0)_{p_n}} \ot e_{\al_n}^{(1)_{p_n}}\ot ... \ot e_{\al_n}^{{(p_n+1)}_{p_n}}\right)^{\Z_m}
$$
where
$$
m={n(p_1+\ldots +p_n +2n)}
$$

A remarkable an almost obvious fact is that the linear subspace
$$
u(Cyc(\widehat{W}_N))\subset Cyc(\widehat{\bW}_N)
$$
 is closed with respect to all strong homotopy involutive Lie bialgebra operations and hence is itself a $\HoLoB$-algebra. In this way we induce a $\HoLoB$-algebra structure

$$
\rho\left(
\Ba{c}\resizebox{16mm}{!}{\xy
(-9,-6)*{};
(0,0)*+{a}*\cir{}
**\dir{-};
(-5,-6)*{};
(0,0)*+{a}*\cir{}
**\dir{-};
(9,-6)*{};
(0,0)*+{a}*\cir{}
**\dir{-};
(5,-6)*{};
(0,0)*+{a}*\cir{}
**\dir{-};
(0,-6)*{\ldots};
(-10,-8)*{_1};
(-6,-8)*{_2};
(10,-8)*{_n};
(-9,6)*{};
(0,0)*+{a}*\cir{}
**\dir{-};
(-5,6)*{};
(0,0)*+{a}*\cir{}
**\dir{-};
(9,6)*{};
(0,0)*+{a}*\cir{}
**\dir{-};
(5,6)*{};
(0,0)*+{a}*\cir{}
**\dir{-};
(0,6)*{\ldots};
(-10,8)*{_1};
(-6,8)*{_2};
(10,8)*{_m};
\endxy}\Ea
\right): \ot^n (Cyc(\widehat{W}_N)) \lon \ot^n(Cyc(\widehat{W}_N))[3-(m+n+2a)]
$$
on $Cyc(\widehat{W}_N)$. It is immediate to see that
 all these operations
have degree $-1$ with respect to the weight-grading by the lengths of cyclic words.
On the linear subspace
 $$
 Cyc(W_N)\subset Cyc(\widehat{W}_N)
 $$
this $\HoLoB$-algebra structure reduces precisely to
Schedler's necklace
Lie bialgebra structure corresponding to the quiver (\ref{1: quiver for W_N}).

\bip


{\large
\section{\bf A geometric application: a new family of $\HoLoB$-operations in string topology}
}

\subsection{Poincare duality algebras} Let $A$ be a finite-dimensional differential non-negatively graded unital commutative algebra. We denote  the (degree zero) multiplication map in $A$ by center dot,
$$
\Ba{ccc}
 A\odot A & \lon & A[0]\\
  a\ot b  & \lon & a\cdot b
\Ea
$$
A dg algebra $A$ as above  is called a {\em dg Poincare duality algebra  of degree $n$}\, if it comes  equipped with a degree $-n$
 {\em orientation map}\,
 $$\fo: A\rar \K[-n]
 $$
 such that the induced graded symmetric paring, or {\it scalar product},
 \Beq\label{6: scalar product in A}
\Ba{rccc}
\langle\ ,\ \rangle: & A\odot A & \lon & \K[-n]\\
                     &a\odot b     & \lon & \fo(a\cdot b)
\Ea
\Eeq
is non-degenerate, and, moreover, $\fo(da)=0$ for any $a\in A$. This scalar product induces a canonical isomorphism
$$
\imath: A\lon A^*[-n]:=\Hom(A,\K)[-n]
$$
such that for any $a,b\in A$ a one has $\imath(a)(b)=\fo(a\cdot b)$. This isomorphism combined with the dualization of the multiplication map induces in turn a degree $n$ diagonal on $A$,
$$
\Ba{rccc}
\Delta: & A & \lon & A\odot A [n]\\
    & a & \lon & \Delta(a)=: \sum a' \ot a''
\Ea
$$
which satisfies
$$
\Delta (a\cdot b)= \sum (a\cdot b') \ot b''=\sum (-1)^{|a||b'|} b'\ot (a\cdot b'') =\sum \sum (-1)^{|a''||b|} (a'\cdot b) \ot a''= \sum  a'\ot (a''\cdot b)
$$
for any $a,b\in A$, and hence makes $A$ into a {\em Frobenius algebra}.
 In particular, we have the equalities (presented up to the obvious Koszul sign as in the formulae just above)
\Beq\label{6: o(abc)}
\sum \fo(a\cdot b\cdot c') c''= \pm \fo(a\cdot b\cdot c'') c'= \pm \fo(a'\cdot b\cdot c) a''=  \pm \fo(a''\cdot b\cdot c) a'=\pm \fo(a\cdot b'\cdot c) b''=\pm \fo(a\cdot b''\cdot c) b'.
\Eeq
for any $a,b,c\in A$, which we shall use later.

\sip

We assume from now on that $A$ is connected and augmented,
$$
A=\K\oplus \bar{A}
$$
 with $\bar{A}$ being positively graded.
 In this case the composite map
\Beq\label{6: diamondsuit}
\diamondsuit: A\stackrel{\Delta}{\lon} A\ot A[n] \stackrel{\cdot}{\lon} A[n]
\Eeq
vanishes when restricted to the subspace $\bar{A}$.

\sip

 As $A$ is finite-dimensional, all the above properties (with appropriate degree changes) hold true for its dual $A^*$, i.e. $A^*$ is also Frobenius algebra but with the diagonal in degree zero and the multiplication map in degree $+n$.

\sip

Given a compact $n$-dimensional manifold $M$, a {\em Poincare model for $M$}\, is, by definition,   a dg Poincare duality algebra $A$ of degree $n$ which is quasi-isomorphic (as an $\assin$ algebra) to the de Rham algebra $\Omega^\bu_M$ of differential forms on $M$. It was proven in \cite{LS} that every connected and simply connected manifold $M$ admits a Poincare model.

\subsection{Cyclic Hochschild complexes and equivariant (co)homology of free loop spaces} Recall (see \S {\ref{2: subsection on cyclic words}} or Corollary {\ref{4: Corollary about HoLoBd}}) that any linear map of the form
$$
\Theta_2: \ot^2(W[d])_{\Z_2} \lon \K[1+d], \ \ \ d\in \Z,
$$
 makes the (reduced) space of cyclic words in $W$
 $$
Cyc(W)=\oplus_{k\geq 0} (\ot^k W)_{\Z_k}\ \  \text{or}\ \  \overline{Cyc}(W)=\oplus_{k\geq 1} (\ot^k W)_{\Z_k},
$$
into a
$\LoBd$-algebra. Let us apply this observation to an arbitrary Poincare duality algebra $A$ in degree $n$. Note that the scalar product (\ref{6: scalar product in A})
induces a non-degenerate graded symmetric linear map,
$$
\Ba{rccc}
\Theta_2': & \ot^2\left(A[n]\right)_{\Z_2} & \lon & \K[n]\\
           & \fs^{n} a \ot \fs^{n} b & \lon & \fs^{n} o(a\cdot b)
\Ea
$$
Indeed, its graded commutativity (i.e.\ invariance under $\Z_2\simeq \bS_2$) can be checked directly,
\Beqrn
\Theta'_2(\fs^{n} b \ot \fs^{n} a) &=& (-1)^{|a||b|} \fs^{n} o(a\cdot b) \\
&=& (-1)^{(|\fs^n a|+n)(|\fs^n b|+n)} \fs^{n} o(a\cdot b) \\
&=& (-1)^{n^2 + n(|\fs^{n}\al|+  |\fs^{n} b|) +  |\fs^{n}\al|\cdot |\fs^{n}\be|} \Theta_2(\fs^{n}a \ot \fs^{n}b)\\
&=& (-1)^{n^2 - n^2 +  |\fs^{n}a|\cdot |\fs^{n}b|} \Theta_2(\fs^{n}a \ot \fs^{n}b)\\
&=& (-1)^{|\fs^{n}\al|\cdot |\fs^{n}b|} \Theta'_2(\fs^{n}a \ot \fs^{n}b)\\
\Eeqrn
for any $a,b\in A$. Here we used the fact that both sides are zero unless $|\fs^{n}a|+  |\fs^{n}b|=-n$. Setting $d:=n-1$ we rewrite the map $\Theta'_2$ as
$$
\Theta'_2: \ot^2\left((A[1])[d]\right) \lon \K[d+1]
$$
and conclude that the spaces of cyclic words $Cyc(A[1])$ and (hence) $\overline{Cyc}(\bar{A}[1])$ come equipped with a $\LoB_{n-1}$-structure, that is, with compatible Lie brackets $\{\ ,\ \}$ and coproduct $\Delta$, both operations having degree $1-d=2-n$.

\sip

The scalar product (\ref{6: scalar product in A})
induces also a non-degenerate linear map,
$$
\Ba{rccc}
\Theta''_2: & \odot^2\left(A[2]\right)\equiv  \ot^2\left((A[n-1])[3-n]\right)_{\Z_2} &\lon &  \K[4-n]=\K[1+(3-n)]\\
 & \imath^{2}(a)\ot \imath^{2}(b)  & \lon & \fs^4 \fo(a\cdot b)
\Ea
$$
implying ---  using the above general construction for $d=3-n$ ---  that the vector spaces
$$
Cyc(A[n-1])\simeq  Cyc(A^*[-1]) \ \ \ \text{and}\ \ \   \overline{Cyc}(\bar{A}^*[-1])
$$
come equipped with an induced $\LoB_{3-n}$-algebra structures, with both basic operations, $\{\ ,\ \}$ and $\Delta$, having degree $1-(3-n)=n-2$. Here we used the canonical isomorphism $A[n]=A^*$.
It is precisely the {\em dual} (in the obvious sense)) of the $\LoB_{n-2}$-algebra structure on
$Cyc(A[1])$ and, respectively,  $\overline{Cyc}(\bar{A}[1])$ discussed just above.

\sip



\subsubsection{\bf Hochschild differential}
The dual of the multiplication map $A\ot A \rar A$
$$
\Ba{rccc}
\triangle: & A^* & \lon & A^*\ot A^* \\
           & a & \lon & \sum a'\ot a''
\Ea
$$
makes the vector space $Cyc(A^*[-1])$ into a complex. If we represents pictorially an element $(\fs^{-1} a_1\ot \fs^{-1} a_2\ot \ldots\ot\fs^{-1}a_i\ot\ldots \ot \fs^{-1} a_n)\in Cyc(A^*[-1])$ as a circle-like  graph with $n$ vertices,
$$
\Ba{c}\resizebox{26mm}{!}{
\begin{xy}
*\xycircle<30pt,30pt>{};
%
 <-30pt,0mm>*{\bu};
  <-25pt,15pt>*{\bullet};
<-25pt,15pt>*{\bullet};
<30pt,0mm>*{\bullet};
  <-25pt,-17pt>*{\bullet};
<25pt,-17pt>*{\bullet};
  <25pt,15pt>*{\bullet};
<25pt,15pt>*{\bullet};
  <-11pt,27pt>*{\bullet};
<11pt,27pt>*{\bullet};
  <-11pt,-28pt>*{\bullet};
<11pt,-28pt>*{\bullet};
<-37pt,0pt>*{^{a_1}};
<-32pt,15pt>*{^{a_2}};
<-32pt,-17pt>*{^{a_n}};
<13mm,5mm>*{^{a_{i-1}}};
<14mm,-0.5mm>*{^{a_{i}}};
<13mm,-6.5mm>*{^{a_{i+1}}};
\end{xy}
}\Ea
$$
with $n$ vertices decorated by elements of $A^*[-1]$ and $n$ edges, then the
 Hochschild differential  $d_H$ in $Cyc(A^*[-1])$ is defined as an alternating sum over splitting each vertex into two new vertices as explained pictorially below,
 $$
d_H:\ \Ba{c}\resizebox{26mm}{!}{
\begin{xy}
*\xycircle<30pt,30pt>{};
%
 <-30pt,0mm>*{\bu};
  <-25pt,15pt>*{\bullet};
<-25pt,15pt>*{\bullet};
<30pt,0mm>*{\bullet};
  <-25pt,-17pt>*{\bullet};
<25pt,-17pt>*{\bullet};
  <25pt,15pt>*{\bullet};
<25pt,15pt>*{\bullet};
  <-11pt,27pt>*{\bullet};
<11pt,27pt>*{\bullet};
  <-11pt,-28pt>*{\bullet};
<11pt,-28pt>*{\bullet};
<-37pt,0pt>*{^{a_1}};
<-32pt,15pt>*{^{a_2}};
<-32pt,-17pt>*{^{a_n}};
<13mm,5mm>*{^{a_{i-1}}};
<14mm,-0.5mm>*{^{a_{i}}};
<13mm,-6.5mm>*{^{a_{i+1}}};
\end{xy}
}
\Ea \lon
\sum_{i=1}^n (-1)^{i}\
\Ba{c}\resizebox{26mm}{!}{
\begin{xy}
*\xycircle<30pt,30pt>{};
%
 <-30pt,0mm>*{\bu};
  <-25pt,15pt>*{\bullet};
<-25pt,15pt>*{\bullet};
<30pt,0mm>*{\bullet};
  <-25pt,-17pt>*{\bullet};
<25pt,-17pt>*{\bullet};
  <25pt,15pt>*{\bullet};
<25pt,15pt>*{\bullet};
  <-11pt,27pt>*{\bullet};
<11pt,27pt>*{\bullet};
  <-11pt,-28pt>*{\bullet};
<11pt,-28pt>*{\bullet};
<-37pt,0pt>*{^{a_1}};
<-32pt,15pt>*{^{a_2}};
<-32pt,-17pt>*{^{a_n}};
<13mm,5mm>*{^{a_{i-1}}};
<14mm,-0.5mm>*{^{a'_{i}}};
<13mm,-6.5mm>*{^{a''_{i}}};
<8mm,-10.5mm>*{^{a_{i+1}}};
\end{xy}
}\Ea
$$

Similarly, the reduced diagonal
$$
\bar{\triangle}: \bar{A}^* \lon A^*  \stackrel{\triangle}{\lon}  A^*\ot A^* \lon \bar{A}^*\ot \bar{A}^*
$$
makes the grade vector space $(\overline{Cyc}(\bar{A}^*[-1]))$ into a complex which sometimes has a nice geometric interpretation (see the next subsection).

\sip

It was shown in \cite{CEG} (see also \cite{CFL, NW} for alternative proofs) that the above mentioned
$\LoB_{3-n}$-algebra structure on $(\overline{Cyc}(\bar{A}^*[-1]))$ is respected by the differential $d_H$ and hence survives at the cohomology level; we conclude that {\em $H^\bu((\overline{Cyc}(\bar{A}^*[-1])), d_H)$ comes equipped with compatible Lie bracket and Lie cobracket of degree $n-2$ for any Poincare duality algebra $A$ in degree $n$}. Similarly, one gets a canonical $\LoB_{n-1}$ algebra structure on $H^\bu((\overline{Cyc}(\bar{A}[1])), d_H^*)$ which has been much studied in \cite{NW} (with many nice explicit formulae obtained).

\subsubsection{\bf $\LoB$ structures on equivariant cohomology of free loop spaces}

Let $LM$ stand for the space of free loops in a closed oriented $n$-dimensional manifold $M$, and
$\bar{C}^\bu_{S^1}(LM)$ (resp., $\bar{C}_\bu^{S^1}(LM)$) for the reduced equivariant cochain  (resp., chain) complex of $LM$ with respect to the obvious $S^1$ action on $LM$; note that we assume that
 $\bar{C}_\bu^{S^1}(LM)$ is the genuine linear dual of $\bar{C}^\bu_{S^1}(LM)$ as a $\Z$-graded complex, i.e. \ it is negatively graded and hence has the boundary differential  of cohomological degree $+1$.
 Let $\bar{H}^\bu_{S^1}(LM)$ (resp.,  $\bar{H}_\bu^{S^1}(LM)$)  stand for its cohomology.

\sip

Assume  $A$ is a Poincare model of $M$, then there are  canonical morphisms of cohomology groups
$$
\bar{H}_\bu^{S^1}(LM) \lon H^{\bu}( \overline{Cyc}(\bar{A}^*[-1]), d_H), \ \ \
H^{\bu}( \overline{Cyc}(\bar{A}[1]), d_H^*) \lon
\bar{H}^\bu_{S^1}(LM)
$$
which are isomorphisms if $M$ is a closed connected and  simply connected manifold. In this case
both the equivariant (co)homology groups
$H_\bu^{S^1}(LM)$ and $H^\bu_{S^1}(LM)$ come equipped with induced involutive Lie bialgebra structure which has been first discovered --- in a very nice geometric way --- by M.\ Chas and D.\ Sullivan in \cite{CS}.

\subsection{A new class of $\HoLoBd$-algebra structures from Poincare duality algebras and props of ribbon hypergraphs} Our main result in this paper is the observation that for any  a graded vector space $W$ equipped with a
 linear map
$$
\Theta_p:  (W[d])^{\ot p}_{\Z_p} \rar \K[1+d]
$$
for some $d\in \Z$ and $p\geq 2$, the associated  vector spaces of ``cyclic words in $W$",
$$
Cyc(W):=\bigoplus_{n\geq 0} (W^{\ot n})^{\Z_n}, \ \ \ \overline{Cyc}(W):=\bigoplus_{n\geq 1} (W^{\ot n})^{\Z_n}
$$
come equipped canonically  with a $\HoLoBd$-algebra structure given by explicit formulae
(\ref{4: rho-diamond}). Given a Poincare duality algebra $A$ in degree $n$ and the associated linear maps
\Beq\label{6: fo_p}
\Ba{rccc}
\fo_p: & \odot^p A & \lon & \K[-n]\\
       & \al_1\ot \cdots \ot \al_p & \lon& \fo( \al_1\cdot \cdots \cdot \al_p),
\Ea
\Eeq
can we use the above construction to induce a $\HoLoBd$-algebra structure on $Cyc(A[1])$ or $Cyc(A^*[-1])$ for suitable $d,p\in \Z$?

\sip

{\sc Case 1}. Consider first the case $W=A[1]$. Out of the map $\fo_p$ we would like to built a map
$$
\Theta_p:  (W[d])^{\ot p}_{\Z_p} \rar \K[1+d].
$$
We can rewrite (\ref{6: fo_p}) as
$$
\fo_p: \odot^p (W[-1]) \lon \K[-n]
$$
or equivalently as
$$
\fo_p:\ot^p (W[-1 +m]) \lon \K[-n +mp], \ \ \ \ \  \forall  m\in \Z.
$$
The equations
$$
d=-1+m, \ \ 1+d=-n+mp
$$
imply
$$
n=m(p-1)
$$
Hence the only solution which applies to every Poincare duality algebra (even to the one whose degree is a prime number) is $p=2$ in which case $n=m$, $d=n-1$; hence in this case we recover the Chas-Sullivan $\LoB_{n-1}$-algebra structure on $Cyc(A[1])$ discussed above.

\sip

{\sc Case 1}. Consider next the case $W={A}[n-1]\simeq {A}^*[-1]$. Then we can rewrite (\ref{6: fo_p}) as
$$
\fo_p: \odot^p (W[1-n]) \lon \K[-n]
$$
or equivalently as
$$
\fo_p:\ot^p (W[1-n +m]) \lon \K[-n +mp], \ \ \ \ \ \forall  m\in \Z.
$$
The equations
$$
d=1-n+m, \ \ 1+d=-n+mp
$$
imply
$$
2=m(p-1)
$$
Here we get precisely two solutions. The first solution is $p=2$ implying  $m=2$, $d=3-n$; hence in this case we recover the Chas-Sullivan $\LoB_{3-n}$-algebra structure on $Cyc(A^*[-1])$ discussed above. The second solution is $p=3$, $m=1$, $d=2-n$ which is studied in the following

\subsubsection{\bf Proposition}
 {\em  Given any dg connected Poincare duality algebra $A$ of degree $n$. The graded vector space
$
\overline{Cyc}(\bar{A}^*[-1])
$
comes equipped canonically with a  $\HoLoB_{2-n}$-algebra structure  induced by the map
\Beq\label{6: fo_3}
\Ba{rccc}
\fo_3: & \odot^3 A & \lon & \K[-n]\\
       & a\odot b \odot c & \lon& \fo( a\cdot b \cdot c),
\Ea
\Eeq
whose  only possibly non-trivial operations are controlled by the following four generators of $\HoLoB_{2-n}$,
$$
\Ba{c}\resizebox{8mm}{!}{\xy
(-5,-6)*{};
(0,0)*+{_0}*\cir{}
**\dir{-};
(0,-6)*{};
(0,0)*+{_0}*\cir{}
**\dir{-};
(5,-6)*{};
(0,0)*+{_0}*\cir{}
**\dir{-};
(0,6)*{};
(0,0)*+{_0}*\cir{}
**\dir{-};
(0,8)*{_1};
(-5,-8)*{_1};
(0,-8)*{_2};
(5,-8)*{_3};
\endxy}\Ea\ \ , \ \
\Ba{c}\resizebox{8mm}{!}{\xy
(-5,6)*{};
(0,0)*+{_0}*\cir{}
**\dir{-};
(0,6)*{};
(0,0)*+{_0}*\cir{}
**\dir{-};
(5,6)*{};
(0,0)*+{_0}*\cir{}
**\dir{-};
(0,-6)*{};
(0,0)*+{_0}*\cir{}
**\dir{-};
(0,-8)*{_1};
(-5,8)*{_1};
(0,8)*{_2};
(5,8)*{_3};
\endxy}\Ea\ \ , \ \
\Ba{c}\resizebox{7mm}{!}{\xy
(-5,-6)*{};
(0,0)*+{_0}*\cir{}
**\dir{-};
(5,6)*{};
(0,0)*+{_0}*\cir{}
**\dir{-};
(5,-6)*{};
(0,0)*+{_0}*\cir{}
**\dir{-};
(-5,6)*{};
(0,0)*+{_0}*\cir{}
**\dir{-};
(-5,8)*{_1};
(5,8)*{_2};
(-5,-8)*{_1};
(5,-8)*{_2};
\endxy}\Ea\ \ , \ \
\Ba{c}\resizebox{2.5mm}{!}{\xy
(0,-6)*{};
(0,0)*+{_1}*\cir{}
**\dir{-};
(0,6)*{};
(0,0)*+{_1}*\cir{}
**\dir{-};
(0,8)*{_1};
(0,-8)*{_1};
\endxy}\Ea
$$
which in turn are controlled via the canonical morphism (\ref{4: rho-diamond}) by the following (respectively) ribbon hypergraphs from  $\RH_{2-n}$,
$$
\xy
(0,4)*{\ast}="1";
(-5,-4)*{\circ}="2";
(0,-4)*{\circ}="3";
(5,-4)*{\circ}="4";
\ar @{-} "1";"2" <0pt>
\ar @{-} "1";"3" <0pt>
\ar @{-} "1";"4" <0pt>
\endxy\ \ ,  \
\ \ \ \xy
(0,5)*{\ast}="1";
(0,-4)*{\circ}="3";
"1";"3" **\crv{(4,0) & (4,1)};
"1";"3" **\crv{(-4,0) & (-4,-1)};
\ar @{-} "1";"3" <0pt>
\endxy\ \ ,  \ \
\xy
(0,5)*{\ast}="1";
(-5,-4)*{\circ}="3";
(5,-4)*{\circ}="4";
"1";"3" **\crv{(-3,5) & (5,4)};
"1";"3" **\crv{(-5,2) & (-5,2)};
\ar @{-} "1";"4" <0pt>
\endxy\ \ ,
\ \ \ \xy
(0,5)*{\ast}="1";
(0,-4)*{\circ}="3";
"1";"3" **\crv{(-5,2) & (5,2)};
"1";"3" **\crv{(5,2) & (-5,2)};
"1";"3" **\crv{(-7,7) & (-7,-7)};
\endxy
%
%
%
%
$$
with vertices and boundaries appropriately (skew)symmetrized (so that their labelling is omitted).

\sip

Moreover, these four operations respect the Hochschild differential $d_H$ in  $\overline{Cyc}(\bar{A}^*[-1])$ and hence induce a  $\HoLoB_{2-n}$-algebra structure on the cohomology
$H^\bu(\overline{Cyc}(\bar{A}^*[-1]), d_H)$.
}

\begin{proof} We have to prove only the last statement in the above Proposition as the rest is simply a direct application of Corollary {\label{4: Corollary about HoLoBd}} to the case $p=3$ and $W=\bar{A}^*[-1]$.

\sip

Let us first study an action of the graph $\Ga_1:=\Ba{c}\resizebox{9mm}{!}{\xy
(0,4)*{\ast}="1";
(-5,-4)*{\circ}="2";
(0,-4)*{\circ}="3";
(5,-4)*{\circ}="4";
\ar @{-} "1";"2" <0pt>
\ar @{-} "1";"3" <0pt>
\ar @{-} "1";"4" <0pt>
\endxy}\Ea
$ on $\overline{Cyc}(\bar{A}^*[-1])$ via the canonical representation $\rho_{\Theta_3}$ introduced  in the proof of Theorem {\ref{3: Prop on repr of Rgra in CycW}} for the particular case when $W=\bar{A}^*[-1]$. Up to the standard Koszul signs this action,
$$
\rho_{\Theta_3}(\Ga_1): \overline{Cyc}(\bar{A}^*[-1])  \bigotimes \overline{Cyc}(\bar{A}^*[-1]) \bigotimes \overline{Cyc}(\bar{A}^*[-1]) \lon \overline{Cyc}(\bar{A}^*[-1])[2n-3]
$$
 is given by the following formula  (cf. \S {\ref{2: subsection on cyclic words}})
$$
\rho_{\Theta_3}(\Ga_1): \left(a_1\ot...\ot a_p\right)_{\Z_p} \bigotimes
\left(b_1\ot...\ot b_q\right)_{\Z_q}
\bigotimes
\left(c_1\ot...\ot c_r\right)_{\Z_r}
\lon \sum_{i=1}^p\sum_{j=1}^p\sum_{k=1}^r \hspace{50mm}
$$
$$
\hspace{10mm}
 \pm
\fo\left(\fs^{1-n}a_i \cdot \fs^{1-n}b_j\cdot \fs^{1-n}c_k\right) \left(a_{i+1}\ot ...\ot  a_{i-1}\ot c_{k+1}\ot ... \ot c_{k-1}\ot b_{j+1}\ot ... \ot b_{j-1}\right)^{\Z_{N}}
$$
where we used the isomorphism ${A}^*[-1]\simeq {A}[n-1]$ and set $N=p+q+r-3$.
The Hochschild differential splits inputs as follows
$$
a_i \rar \sum a_i'\ot a_i'', \ \ \ b_j\rar \sum {b_j'\ot b_j''}, \ \ \ c_k\rar \sum c_k'\ot c_k''.
$$
Hence the obstructions for the above map to respect the Hochschild differential $d_H$ are given by the summands of the form,
$$
\fo\left(\fs^{1-n}a_i' \cdot \fs^{1-n}b_j\cdot \fs^{1-n}c_k\right)(...\ot b_{j-1}\ot a_i''\ot a_{i+1}\ot...) + \fo\left(\fs^{1-n}a_i'' \cdot \fs^{1-n}b_j\cdot \fs^{1-n}c_k\right)(...\ot a_{i-1}\ot a_i''\ot c_{k+1}\ot ...)+
$$
$$
\pm\fo\left(\fs^{1-n}a_i \cdot \fs^{1-n}b_j'\cdot \fs^{1-n}c_k\right)(...\ot c_{k-1}\ot b_j''\ot b_{j+1}\ot...) + \fo\left(\fs^{1-n}a_i \cdot \fs^{1-n}b_j''\cdot \fs^{1-n}c_k\right)(...\ot b_{j-1}\ot b_j'\ot a_{i+1}\ot ...)
$$
$$
\pm\fo\left(\fs^{1-n}a_i \cdot \fs^{1-n}b_j\cdot \fs^{1-n}c_k'\right)(...\ot c_{k-1}\ot c_k''\ot b_{j+1}\ot...) + \fo\left(\fs^{1-n}a_i \cdot \fs^{1-n}b_j\cdot \fs^{1-n}c_k''\right)(...\ot a_{i-1}\ot c_k'\ot c_{k+1}\ot ...)
$$
in which precisely one of the ``dashed" elements lands inside the orientation map $\fo$.
Thanks to identity ({\ref{6: o(abc)}}), all these terms come in pairs which cancel each other (for example, the first term in the first row cancels the second term in the second row and so on).

\sip

Consider next the operation on $\overline{Cyc}(\bar{A}^*[-1])$  controlled by the graph  $\Ga_2=\xy
(0,5)*{\ast}="1";
(0,-4)*{\circ}="3";
"1";"3" **\crv{(4,0) & (4,1)};
"1";"3" **\crv{(-4,0) & (-4,-1)};
\ar @{-} "1";"3" <0pt>
\endxy$. The obstructions to this operations under the canonical morphism $\rho_{\Theta_3}(\Ga_2)$ to respect the Hochschild differential $d_H$ are given by the terms of two types,
\Bi
\item[(i)]  one type corresponds to summands in which precisely one of the ``dashed" elements lands inside the orientation map as above; these terms cancel each other due to the same identities ({\ref{6: o(abc)}}) as above;

\item[(ii)]  the second type corresponds to summands in which {\em both}\, ``dashed" elements lands inside the orientation map, i.e.\ to summands  which contain a numerical factor of the form
    $$
    \fo\left(\fs^{1-n}a_i' \cdot \fs^{1-n}a_i''\cdot \fs^{1-n}b\right)
    $$
\Ei
 All these summands  in (ii) vanish dues to the fact that the map (\ref{6: diamondsuit}) vanishes on elements of $\overline{Cyc}(\bar{A}^*[-1])$. The obstructions for the two remaining operations,
$$
\rho_{\theta_3}\left( \xy
(0,5)*{\ast}="1";
(0,-4)*{\circ}="3";
"1";"3" **\crv{(-5,2) & (5,2)};
"1";"3" **\crv{(5,2) & (-5,2)};
"1";"3" **\crv{(-7,7) & (-7,-7)};
\endxy\right) \ \ \text{and}\ \
\rho_{\Theta_3}\left(\xy
(0,5)*{\ast}="1";
(-5,-4)*{\circ}="3";
(5,-4)*{\circ}="4";
"1";"3" **\crv{(-3,5) & (5,4)};
"1";"3" **\crv{(-5,2) & (-5,2)};
\ar @{-} "1";"4" <0pt>
\endxy
\right)
$$
 to respect $d_H$  are  of the two types discussed above; they vanish for the same reasons.

\end{proof}

\subsubsection{\bf Corollary}{\it For any closed connected and simply connected  $n$-dimensional manifold $M$ the reduced equivariant homology
$\bar{H}_\bu^{S^1}(LM)$ of the associated free loop space comes equipped with a $\HoLoB_{2-n}$-algebra structure controlled by the above four operations.}

\sip

It is easy to see that these four new $\HoLoB_{2-n}$ operations on $\bar{H}_\bu^{S^1}(LM)$ are non-trivial for $M=\CP^n$ with $n\geq 3$.

\def\cprime{$'$}

\end{document}